\documentclass[11pt]{amsart}
\usepackage{amsfonts,amssymb,amscd}
\usepackage[all]{xy}

\newtheorem{lemm}{Lemma}[section]
\newtheorem{theo}[lemm]{Theorem}
\newtheorem{coro}[lemm]{Corollary}
\newtheorem{prop}[lemm]{Proposition}

\def\a{\alpha}
\def\JO{J_3(\mathbb{O})}

\def\OP2{\mathbb{OP}^2}
\def\HP2{\mathbb{HP}^2}
\def\bOP2{\bar{\mathbb{OP}}^2}
\def\CC{{\mathbb C}}
\def\HH{{\mathbb H}}
\def\OO{{\mathbb O}}\def\PP{{\mathbb P}}
\def\QQ{{\mathbb Q}}
\def\ZZ{{\mathbb Z}}
\def\ra{\rightarrow}
\def\cO{\mathcal{O}}\def\cD{\mathcal{D}}\def\cA{\mathcal{A}}

\def\cN{\mathcal{N}}\def\cC{\mathcal{C}}

\def\fe{\mathfrak{e}}\def\ff{\mathfrak{f}}
\def\cQ{\mathcal{Q}}
\def\fsl{\mathfrak{sl}}
\def\cS{\mathcal{S}}\def\cD{\mathcal{D}}
\def\cE{\mathcal{E}}\def\cK{\mathcal{K}}\def\cF{\mathcal{F}}
\def\lra{\longrightarrow}

\begin{document}

\title{On cubic hypersurfaces \linebreak of dimension seven and eight}
\author[A. Iliev]{Atanas Iliev}
\address{Department of Mathematics,
Seoul National University,
Seoul 151-747, Korea}
\email{{\tt ailiev2001@yahoo.com}}
\author[L. Manivel]{Laurent Manivel}
\address{Institut Fourier,  
Universit\'e de Grenoble et CNRS,
BP 74, 38402 Saint-Martin d'H\`eres, France}
\email{{\tt Laurent.Manivel@ujf-grenoble.fr}}

\begin{abstract}
Cubic sevenfolds are examples of Fano manifolds of Calabi-Yau type. 
We study them in relation with the Cartan cubic, the $E_6$-invariant
cubic in $\PP^{26}$. We show that a generic cubic sevenfold $X$ can be 
described as a linear section of the Cartan cubic, in finitely
many ways. To each such ``Cartan representation''
we associate a rank nine vector bundle on $X$ with very special
cohomological properties. In particular it allows to define 
auto-equivalences of the non-commutative Calabi-Yau threefold 
associated to $X$ by Kuznetsov. Finally we show that the generic
eight dimensional section of the Cartan cubic is rational. 
\end{abstract}

\maketitle

\section{Introduction}

Manifolds of Calabi-Yau type were defined in \cite{im-fcy} as compact complex
manifolds of odd dimension whose middle dimensional Hodge structure is 
similar to that of a Calabi-Yau threefold. Certain manifolds of Calabi-Yau 
type were used in order to construct mirrors of rigid Calabi-Yau threefolds, 
and in several respects they behave very much like Calabi-Yau threefolds.
This applies in particular to the cubic sevenfold, 
which happens to be the mirror of the so-called 
$Z$-variety -- a rigid Calabi-Yau threefold 
obtained as a desingularization of the quotient of a product of three Fermat 
cubic curves by the group ${\ZZ}_3 \times {\ZZ_3}$, 
see e.g. \cite{sch}, \cite{cdp}, \cite{bbvw}.

In some way, manifolds of Calabi-Yau type come just after manifolds of K3 type:
manifolds of even dimension $2n$, whose middle dimensional Hodge structure is 
similar to that of a K3 surface. For $n \ge 2$, projective $2n$-folds of K3 type 
seem to be scarce, but they are specially interesting in the fact that if one 
has a well-behaved family of $(n-1)$-cycles on the variety, 
then the Abel-Jacobi map produces a closed holomorphic 
two-form on the base of the family. 
This could be a method to construct new complex symplectic 
manifolds. Indeed, examples suggest that these holomorphic forms tend to be as 
non-degenerate as they can, see e.g. \cite{kmm}. 
This is what we will check for the cubic sevenfold
and its family of planes: since we start with a manifold of Calabi-Yau type,
the Abel-Jacobi map will produce a holomorphic three-form on the corresponding 
Fano variety, and we show that it is generically non-degenerate (Proposition 
\ref{nondeg3form}), in a sense that we will make precise. 

Then we will concentrate, for the remaining of the paper, on a phenomenon
which was first observed in \cite{im-fcy}, and that we will study here in 
much more depth: one can construct manifolds of Calabi-Yau type as linear 
sections of the Cayley plane, and this turns out to be again related to cubic sevenfolds. 
The {\it Cayley plane}, which is homogeneous under $E_6$, 
 is the simplest of the homogeneous spaces of exceptional
type. It has many beautiful geometric properties, in particular it supports 
a plane projective geometry, whose ``lines'' are eight dimensional quadrics. 
Its minimal equivariant embedding is inside a projective space of dimension
$26$, in which its secant variety is an $E_6$-invariant cubic hypersurface that
we call the {\it Cartan cubic}. The Cartan cubic is singular exactly along the 
Cayley plane, so that its general linear section of dimension at most 
eight is smooth. 

We prove in Proposition \ref{finiterep} that a general cubic sevenfold $X$
can always be described as a linear section of the Cartan cubic. Moreover,
there is, up to isomorphism, only finitely many such {\it Cartan representations}. 
Using the plane projective geometry supported by the Cayley plane
we show that such a representation induces on $X$ a rank nine vector bundle 
$\cE_X$ with very nice properties: in particular, it is simple, infinitesimally 
rigid (Theorem \ref{rigidity}) and arithmetically Cohen-Macaulay (Proposition \ref{aCM}). 
Kuznetsov observed in \cite{kuz1}
that the derived category of coherent sheaves on a smooth hypersurface contains,
under some hypothesis on the dimension and the degree, a full subcategory which is 
a Calabi-Yau category, of smaller dimension than the hypersurface. In 
particular  the derived category $D^b(X)$ of our smooth cubic sevenfold contains a 
full subcategory $\cA_X$ which is Calabi-Yau of dimension three (a ``{\it non-commutative
Calabi-Yau threefold}''). It turns out that $\cE_X(-1)$ and $\cE_X(-2)$ both define 
objects in  $\cA_X$, and these objects are spherical. As 
shown by Seidel and Thomas, the corresponding spherical twists define 
auto-equivalences of $\cA_X$ (Corollary \ref{sphtwists}). 
Recall that $D^b(X)$ itself is poor in 
symmetries, $X$ being Fano. Since $\cA_X$ is Calabi-Yau, its structure and
its symmetries are potentially much richer. 

Cartan representations of cubic sevenfolds are very similar to Pfaffian 
representations of cubic fourfolds, which are linear sections of the 
secant cubic to the Grassmannian $G(2,6)$. (There are nevertheless 
important differences, in particular a general cubic fourfold is {\it not} 
Pfaffian.) Starting from a Pfaffian representation, a classical construction
consists in considering the orthogonal linear section of the dual Grassmannian:
this is a K3 surface of genus eight, whose derived category can be embedded 
inside the derived category of the fourfold \cite{kuz2}. We expect something  
similar for a Cartan representation of a general cubic sevenfold $X$. 
The orthogonal linear section of the  dual Cayley plane is a Fano 
sevenfold $Y$ of index three. Using the plane projective geometry supported
by $\OP2$,  we show that $Y$ is birationally equivalent (but not 
isomorphic) to $X$ (Proposition \ref{tto}). 
This is an analogue of the Tregub-Takeuchi birational
equivalence between the general prime Fano threefold of degree fourteen
and the orthogonal cubic threefold, as we explain in a short  Appendix.
We conjecture that an explicit subcategory 
$\cA_Y$ of $D^b(Y)$ should be Calabi-Yau of dimension three, and  
should even be equivalent to $\cA_X$. This would give rise to an 
isomorphism, very similar to the Pfaffian-Calabi-Yau derived equivalence 
of \cite{bc}, between two non-commutative Calabi-Yau threefolds,
associated to two quite different sevenfolds of Calabi-Yau type. 

Pfaffian fourfolds are famous examples of {\it rational} cubic fourfolds.
We conclude the paper by showing (Theorem \ref{rational}) 
that eight dimensional linear sections of
the Cartan cubic give rise to rational cubic eightfolds (but very non
generic). The main ingredient of the proof is the close connection between 
$E_6$ and $Spin_{10}$, and the fact that the spinor varieties of $Spin_{10}$
have four dimensional linear sections which are varieties with one 
apparent double point. The very nice geometric and homological properties 
of Pfaffian fourfolds thus seem to propagate in rather different ways 
to cubics of dimension seven or eight. 

\medskip\noindent {\it Acknowledgments}. 
The second author acknowledges financial
support from the ANR Project VHSMOD09 ANR-09- BLAN-0104-01.


\section{Planes on the cubic sevenfold}\label{planes} 

Here we shall see that on the 8-dimensional family $F = F(X)$ of planes 
on the cubic sevenfold $X$ there exists a generically non-degenerate 
holomorphic 3-form $\omega$. The existence of the form 
$\omega \in H^{3,0}(F(X))$ is based on the fact that the space 
$H^{5,2}(X)$ is one-dimensional, and $\omega$ can be constructed for example by 
following the techniques from \cite{kmm}. 
This is an analog of the well-known observation 
about the cubic fourfold $Y$ and its family of lines $F(Y)$, 
due to Beauville and Donagi: the unique (3,1)-form 
(up to scalars) on the cubic fourfold $Y$ yields a non-degenerate 
holomorphic 2-form on $F(X)$, thus 
establishing that $F(X)$ is a holomorphic symplectic fourfold,
see \cite{bd}.

\subsection{Basics}

The cubic sevenfold and its Abel-Jacobi mapping have been studied
by Albano and Collino \cite{alb-col}. It also appeared in the physics 
literature, since it was observed in \cite{cdp} that it could be 
considered as the mirror of a rigid Calabi-Yau threefold, obtained 
as the quotient of a product of three elliptic curves by some 
threefold symmetry group. 

In \cite{im-fcy} we defined {\it manifolds of Calabi-Yau type} as compact
complex manifolds $X$ of odd dimension $2n+1$ such that:
\begin{enumerate}
\item
The middle dimensional Hodge structure is numerically similar to that of a 
Calabi-Yau threefold, that is 
$$h^{n+2,n-1}(X)=1, \qquad \mathrm{and}\quad h^{n+p+1,n-p}(X)=0 \;\;
\mathrm{for}\;p\ge 2.$$
\item
For any generator $\omega\in H^{n+2,n-1}(X)$, the contraction map 
$$H^1(X,TX)\stackrel{\omega}{\lra} H^{n-1}(X,\Omega_X^{n+1})$$
is an isomorphism. 
\item
The Hodge numbers $h^{k,0}(X)=0$ for $1\le k\le 2n$.
\end{enumerate}

It is a straightforward consequence of the description of the Hodge 
structure of a smooth projective hypersurface in terms of its Jacobian 
ring, that a smooth cubic sevenfold $X$ is a Fano manifold of Calabi-Yau type.
The relevant Hodge numbers are
$$h^{5,2}(X)=1, \qquad  h^{4,3}(X)=84.$$

\subsection{A 3-form on the family of planes in the cubic sevenfold}

Our first interest in cubic sevenfolds stemmed from the fact that 
the equality $h^{5,2}(X)=1$ looks very similar to the equality 
$h^{3,1}(Y)=1$ that holds for a smooth cubic fourfold $Y$. 
As we have already mentioned, it is a famous result of Beauville and Donagi 
\cite{bd} that the Fano fourfold $F(Y)$ 
of lines on the general cubic fourfold $Y$ is a symplectic fourfold: 
it is endowed with a closed holomorphic two-form $\omega\in H^{2,0}(F(Y))$
which is everywhere non-degenerate. Moreover, 
if $\mathcal{L}\subset F(Y)\times Y$ is the  family of lines on $Y$, with its projections
$p$ and $q$ to $F(Y)$ and $Y$, then the natural map 
  $$p_*q^* : H^{3,1}(Y)\ra H^{2,0}(F(Y))$$
is an isomorphism.  

If now $X$ is a general cubic sevenfold in $\PP^8$, then the Fano 
variety $F=F(X)$ of {\it planes} on $X$ is a smooth connected eightfold, 
and with the same notations as above, the map
 $$p_*q^* : H^{5,2}(X)\ra H^{3,0}(F)$$
is an isomorphism \cite{alb-col}. In particular $H^{3,0}(F)=\CC$, so that $F$ is endowed 
with a holomorphic three-form $\omega$, which is closed by the results 
of \cite{km}. Although $F$ is not a symplectic variety
(its canonical bundle is non-trivial),
we have a partially similar result to the case of the cubic fourfold. 

\begin{prop}\label{nondeg3form}
Let $X$ be a general cubic sevenfold, $F$ the eightfold of planes
on $X$, and $\omega$ a generator of  $H^{3,0}(F)$. Then $\omega$
is generically non-degenerate. 
\end{prop}

Before proving this statement, a few words of explanation are in 
order. Indeed, if the concept of non-degeneracy is completely clear 
for two-forms, it is not for three-forms in general. The point here 
is that $F$ has dimension eight, so that $\omega$ is at each point
of $F$ a three-form in eight variables. Now, it is known that the 
space of skew-symmetric three forms in eight variables has a finite
number of $GL_8$-orbits, which were first classified by Gurevich 
\cite{gu}. In particular there 
is an open orbit 
$\cO_{gen}\subset \wedge^3\CC^8$, and $\omega$ can be defined as 
being non-degenerate at some point $P$ if 
$$\omega\in\cO_{gen}\subset \wedge^3(T_PF)^\vee.$$

\proof Following \cite{km}, 
we can describe the three-form in terms of the extension
class $\nu\in Ext^1(\Omega^1_X,\cO_X(-3))$ defining the 
conormal exact sequence. Its square is an element $S^2\nu\in 
Ext^2(\Omega^2_X,\cO_X(-6))$. 
Let $P$ be a general plane in $X$. 
The normal bundle $N_{P/X}$ has rank five and determinant 
$\wedge^5N_{P/X}=\cO_P(3)$. Moreover, the quotient $TX_{|P}\ra N_{P/X}$
can be interpreted as a global section $t$ of $\Omega^1_{X|P}\otimes N_{P/X}$. 
We can thus define the composition 
\begin{eqnarray*}
 &\wedge^3H^0(N_{P/X})\ra H^0(\wedge^3N_{P/X})\stackrel{\wedge^2t}{\ra}
H^0(\wedge^5N_{P/X}\otimes\Omega^2_{X|P})= \\
 &\hspace*{5cm}=H^0(\Omega^2_{X|P}(3))\stackrel{S^2\nu}{\ra}H^2(\omega_P)=\CC.
\end{eqnarray*}
Once we identify $T_PF$ with $H^0(N_{P/X})$, this is precisely the three-form
$\omega$ (up to scalar) and we can make an explicit local computation.

\smallskip
So let $\CC^9=A\oplus B$ be a direct sum decomposition such that the plane $P$ 
is the projectivization of $A$. Let $x_1,x_2,x_3$ be coordinates on $A$ and 
$y_1,\ldots , y_6$ be coordinates on $B$. Let $c(x,y)$ be a cubic polynomial 
defining the hypersurface $X$. The hypothesis that $X$ contains $P$ implies
that $c$ can be written as 
$$c(x,y)=\sum_{i=1}^6y_iq_i(x,y)$$
for some quadratic forms $q_1,\ldots ,q_6$. The normal bundle $N_{P/X}\subset 
B\otimes\cO_P(1)$ can be described as the set of sixtuples $\ell_1(x),\ldots ,
\ell_6(x)$ such that 
$$\sum_{i=1}^6\ell_i(x)q_i(x,0)=0.$$
The space of global sections $H^0(N_{P/X})$ has then dimension eight if and 
only if any cubic form on $A$ can be written in the previous form, that is, 
$\sum_{i=1}^6\ell_i(x)q_i(x,0)$ for some linear forms $\ell_i(x)$. In the sequel 
we make a stronger assumption: we suppose that the quadrics $q_1(x,0)
,\ldots ,q_6(x,0)$ are linearly independent. Otherwise said, they define an 
isomorphism between $B$ and $S^2A$. 

This suggests to change our notation
by identifying $i=1\ldots 6$ with pairs of integers $(j,k)$ from $1$ to $3$, 
and letting $q_{jk}(x,0)=x_jx_k$. Correspondingly, $B$ has a basis $e_{jk}$
and $N_{P/X}\subset B\otimes\cO_P(1)$ is the set of sixtuples $(\ell_{jk}(x))$
such that $\sum_{jk}\ell_{jk}(x)x_jx_k=0$. An easy computation shows that 
$\wedge^5N_{P/X}\subset \wedge^5B\otimes\cO_P(5)\simeq B^\vee(5)$ is the 
copy of $\cO_P(3)$ generated by the global section $n(x)=\sum_{jk}x_jx_ke_{jk}^\vee$
of  $B^\vee(2)$, which never vanishes.

\smallskip
Note that the map $TX_{|P}\ra N_{P/X}$ induces a map from $Ext^1(\Omega^1_X,\cO_X(-3))$
to $Ext^1(N_{P/X}^\vee,\cO_X(-3))$, mapping the class $\nu$ to the class $\nu_P$ 
of the extension 
\begin{equation}\label{ext}
0\ra  N_{P/X}\ra B\otimes\cO_P(1)\ra \cO_P(3)\ra 0.
\end{equation}
This allows to compute our three-form on $H^0(N_{P/X})$ as the composition 
\begin{eqnarray*}
 &\wedge^3H^0(N_{P/X})\ra H^0(\wedge^3N_{P/X})\hspace*{2cm} \\
 &\hspace*{4cm} \simeq H^0(\wedge^2N_{P/X}^\vee\otimes \cO_P(3))
  \stackrel{S^2\nu_P}{\ra}H^2(\omega_P)=\CC.
\end{eqnarray*}
Let us find a Cech representative of $\nu_P$ relative to the canonical affine 
covering of $P$. On $x_i\ne 0$, the section $x_i\otimes e_i^2$ yields a splitting 
of $(\ref{ext})$, which implies that $\nu_P$ is defined by the Cech cocycle 
$$g_{ij}=\frac{x_i\otimes e_{ii}}{x_i^3}-\frac{x_j\otimes e_{jj}}{x_j^3}.$$
The symmetric square of $\nu_P$ in $Ext^2(\wedge^2N_{P/X}^\vee,\cO_X(-6))$ is defined
by the Cech cocycle 
$$g_{ijk}=\frac{e_{ii}\wedge e_{jj}}{x_i^2x_j^2}+\frac{e_{jj}\wedge e_{kk}}{x_j^2x_k^2}
+\frac{e_{kk}\wedge e_{ii}}{x_k^2x_i^2}.$$
Let now $a,b,c$ be three global sections of $N_{P/X}$. For example, we can 
identify $a$ with a traceless matrix $(a_{jk})$ of order three, in such a way
that, as a global section of $B(1)$, 
$$a=\sum_{\sigma\in\cS_3} \sum_{k=1}^3
 \epsilon(\sigma)a_{\sigma(1)k}x_{\sigma(2)}e_{\sigma(3)k}.$$
Then our three form $\psi$ on $H^0(N_{P/X})$ maps $a\wedge b\wedge c$ to
the Cech cocycle of $\omega_P$ defined by $a\wedge b\wedge c\wedge g_{123}$, which 
can be written, when considered as a rational section of $\wedge^5N_{P/X}$,
as a multiple of $n(x)$ by some rational function having poles at most along the three
coordinate hyperplanes. Since $1/x_1x_2x_3$ defines a cocycle 
generating $H^2(\omega_P)$, while any other monomial (with exponents in $\ZZ$) 
is a coboundary, we deduce the identity
$$a\wedge b\wedge c =\psi(a,b,c)x_1x_2x_3 e_{12}\wedge e_{23}\wedge e_{13} +\mathrm{other\;terms}.$$
At that point the simplest thing to do is an explicit computation: we find, in terms of the 
linear forms $\alpha_{jk}$ on $\fsl_3$ taking the values of the $(jk)$-entries (so that 
$ \alpha_{11}+\alpha_{22}+\alpha_{33}=0$):
\begin{eqnarray*}
 &\psi = &\alpha_{12}\wedge\alpha_{31}\wedge\alpha_{23} +\alpha_{13}\wedge\alpha_{21}\wedge\alpha_{32} +
 \alpha_{12}\wedge\alpha_{21}\wedge (\alpha_{22}-\alpha_{11}) +\\
 & & +\alpha_{23}\wedge\alpha_{32}\wedge (\alpha_{33}-\alpha_{22})+
\alpha_{13}\wedge\alpha_{31}\wedge (\alpha_{11}-\alpha_{33}),
 \end{eqnarray*}
in which we easily recognize the unique invariant three-form on $\fsl_3$:
$$\psi(a,b,c)=\mathrm{trace}([a,b]c).$$
Any element in the stabilizer of $\psi$ clearly induces a Lie algebra isomorphism of $\fsl_3$. This stabilizer
is thus locally isomorphic to $SL_3$. In particular it has dimension $8$, hence codimension
$56$ in $GL(\fsl_3)$. But $56$ is also the dimension of $\wedge^3\fsl_3$, which shows that 
the $GL(\fsl_3)$-orbit of the determinant is the open orbit $\cO_{gen}$ of non-degenerate 
three-forms. 
\qed

\medskip\noindent {\it Remark}. Following \cite{sk} (Proposition 10 page 90), 
the complement of $\cO_{gen}$ is defined
by a covariant of degree $16$ on  $\wedge^3\CC^8$. This covariant defines a map $S^{16}\Omega_{F(X)}^3\ra
\omega_{F(X)}^6$. Since $\omega_{F(X)}=\cO_{F(X)}(1)$, the degeneracy locus of the three-form 
 $\omega$ is a sextic hypersurface in $F(X)$.


\section{Cartan representations of cubic sevenfolds}\label{cartan}

By a result of L. Hesse from 1844 the homogeneous form defining a 
smooth plane cubic can be represented in three non-equivalent ways as 
a determinant of a symmetric $3 \times 3$ matrix with linear entries, 
see \cite{he}. 
By a later construction of A. Dixon (1902), the representations 
of the equation of a plane curve as a symmetric determinant 
correspond to its non-effective theta-characteristics, 
see \cite{di} or \S 4 in \cite{beau}.  
By other classical results originating from A. Clebsch (1866), 
a smooth cubic surface $S\subset \PP^3$ 
has exactly $72$ determinantal representations, 
i.e. the cubic form $F(x)$ of $S$ has $72$  non-equivalent representations 
as a determinant of a $3 \times 3$ matrix $M$ with linear entries, 
see \cite{cl} or Cor.6.4 in \cite{beau}.  

What about cubics of higher dimensions? 
The determinantal cubics of dimensions $\ge 3$ are singular,
but one can nevertheless proceed by replacing determinantal 
by Pfaffian representations, i.e. representations as Pfaffians 
of skew-symmetric matrices of order $6$ with linear entries. 
It has been found essentially in \cite{mt}, \cite{imr} and \cite{dr}, 
and restated in this form in \cite{beau}, that the general 
cubic threefold $Y$ has a family of Pfaffian representations
which is 5-dimensional and birational to the intermediate 
Jacobian of $Y$. 
As for cubic fourfolds, it is well-known 
that the general cubic fourfold is not Pfaffian. 
For cubics of higher dimensions the situation is even worse 
-- while the general Pfaffian cubic fourfold is still smooth
(and rational !) the Pfaffian cubics of dimension greater 
than five are singular, and the question is which generalization 
of the Pfaffian, if any, can represent cubics of higher dimensions. 

Looking back, we see that the spaces of symmetric $3 \times 3$ matrices,
matrices of order 3, and skew-symmetric matrices of order 6 
are in fact the first three of a series of four Jordan algebras --
they are the algebras of Hermitian matrices over the composition algebras 
${\bf R}$, ${\bf C}$ and ${\bf H}$. 
Then, writing-down  a cubic form $F$ as a symmetric 
determinant,  as a determinant, or as a Pfaffian is the same as 
giving a 
presentation of the cubic $F = 0$ as a linear section of the 
determinant in the Jordan algebra over ${\bf R}$, ${\bf C}$ or ${\bf H}$. 

What is left is the fourth Jordan algebra -- the space of 
Hermitian  $3 \times 3$ matrices over the Cayley algebra ${\bf O}$ 
of octonions. 
The octonionic determinant is a special $E_6$-invariant 
cubic in $\PP^{26}$ -- the Cartan cubic which is singular 
in codimension 9. 
Therefore we can ask about representations of the 
smooth cubics $X$ as linear sections of the Cartan cubic upto 
dimension 8. We call them the {\it Cartan} representations of $X$. 

As we show below, the general cubic sevenfold is a section of the Cartan 
cubic in a finite number of ways, and by this property it looks 
similar to the cubic curve and the cubic surface. 
But while the numbers of determinantal representations as above 
of cubic curves and of cubic surfaces are known 
since Hesse and Clebsch, the number of Cartan 
representations of the general cubic sevenfold remains an unsolved 
question. 


\subsection{The Cayley plane and the Cartan cubic}

Let ${\bf O}$ denote the normed algebra of real octonions (see e.g. \cite{baez}), and let 
$\OO$ be its complexification. The space
$$\JO = \Bigg\{ x=
\begin{pmatrix} 
r            & w            & \bar{v}  \\ 
\bar{w} & s            &  u            \\ 
v            & \bar{u} &  t            \\ 
\end{pmatrix}
: r,s,t \in {\bf C}, u,v,w\in \OO \Bigg\} \cong {\bf C}^{27}
$$
of $\OO$-Hermitian matrices of order $3$, is traditionally
known as the {\it exceptional 
simple Jordan algebra}, for the Jordan multiplication 
$A \circ B = \frac{1}{2}(AB+BA)$. The automorphism 
group of this Jordan algebra is an algebraic group of
type $F_4$, and preserves the determinant 
$$Det(x) = rst-r|u|^2-s|v|^2-t|w|^2
+2\mathrm{Re}(uvw),$$
where the term $\mathrm{Re}(uvw)$ makes sense uniquely, 
despite the lack of associativity of $\OO$. 
The subgroup of $GL(\JO)$ consisting of automorphisms preserving 
this cubic polynomial is the adjoint group of type $E_6$. 
The Jordan algebra $\JO$ and its dual are the minimal 
(non-trivial and irreducible) representations of this group. 
An equation of the invariant cubic hypersurface was once written
down by Elie Cartan in terms of the configuration of the $27$ lines
on a smooth cubic surface. We will call it the {\it Cartan cubic}.

The action of $E_6$ on the 
projectivization $\PP\JO$ has exactly three orbits: the complement of the 
determinantal hypersurface $(Det=0)$, the regular part of this hypersurface, and its singular
part which is the closed $E_6$-orbit. These three orbits can be considered as 
the (projectivizations of the) sets 
of matrices of rank three, two, and one respectively. 

The closed orbit, corresponding to rank one matrices,
is called the {\it Cayley plane} and denoted $\OP2$ (the reasons for this notation and terminology 
will soon be explained). It can be defined by the quadratic equation 
$$x^2=\mathrm{trace}(x)x, \qquad  x\in \JO.$$ 
It will be useful to notice that a dense open subset can be parametrized explicitly 
as the set of matrices of the form 
\begin{equation}\label{dense}
\begin{pmatrix} 
1           & w            & \bar{v}  \\ 
\bar{w} & |w|^2            &  \bar{w}\bar{v}          \\ 
v            & vw &  |v|^2               \\ 
\end{pmatrix}, 
 \quad v,w\in \OO . 
\end{equation}
The Cayley plane  can also be identified with the 
quotient of $E_6$ by the maximal parabolic subgroup $P_1$ defined 
by the simple root $\a_1$ corresponding to the end of one the 
long arms of the Dynkin diagram (we follow
the notations of \cite{bou}). The semi-simple part 
of this maximal parabolic is isomorphic to ${\rm Spin}_{10}$.
\medskip

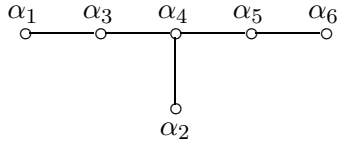
\begin{center}
\setlength{\unitlength}{5mm}
\begin{picture}(20,4)(-5,1)
\multiput(-.3,3.8)(2,0){4}{$\circ$}
\multiput(0,4)(2,0){4}{\line(1,0){1.7}}
\put(3.7,1.8){$\circ$}
\put(7.7,3.8){$\circ$}
\put(-.6,4.4){$\a_1$}
\put(1.4,4.4){$\a_3$}
\put(3.4,4.4){$\a_4$}
\put(5.4,4.4){$\a_5$}
\put(7.4,4.4){$\a_6$}
\put(3.45,1.3){$\a_2$}
\put(3.85,2.1){\line(0,1){1.8}}
\end{picture}
\end{center}

The symmetric node  $\a_6$ of the Dynkin diagram corresponds to the 
dual representation $\JO^\vee$, which is exchanged with $\JO$
by an outer automorphism of the group. In particular the 
orbits inside $\PP\JO$ and $\PP\JO^\vee$ are isomorphic 
(although non-equivariantly), and the 
closed orbit $E_6/P_6$ inside $\PP \JO^\vee$ is again a 
copy of the Cayley plane. We denote it by $\bOP2$
and call it the {\it dual Cayley plane}.

The Cartan cubic inside $\PP \JO^\vee$ can be described 
either as the projective dual variety to the Cayley plane,
or as the secant variety to the dual Cayley plane. As a 
secant variety, it is degenerate, and therefore any
point $p$ of the open orbit in  $\PP \JO^\vee$ defines a
positive dimensional {\it entry locus} in the dual Cayley
plane (the trace on $\bOP2$ of the secants through
$p$, see e.g. \cite{zak} for background and more details
on Severi varieties). 
One can check that these entry-loci are eight
dimensional quadrics, and that their family is parametrized
by the Cayley plane itself. Explicitely, if $x$ is in 
$\OP2$, the orthogonal to its tangent space cuts the 
dual Cayley plane along one of these quadrics $Q_x$. 

Quite remarkably, these quadrics must be considered as 
$\OO$-lines inside the dual Cayley plane, and we get the 
familiar picture of a plane projective geometry: the two
Cayley planes parametrize points and lines (respectively 
lines and points) on them, and the basic axioms of a plane 
geometry hold generically \cite{tits, freud}:
\begin{enumerate}
\item Let $x,x'\in\OP2$ be two distinct points, such that 
the line $xx'$ is not contained in $\OP2$. 
Then $Q_x\cap Q_{x'}$ is a unique reduced point. 
\item Let $y,y'\in\bOP2$ such that the  
line $yy'$ is not contained in $\bOP2$. Then there
is a unique $x$ in $\OP2$ such that the $\OO$-line
$Q_x$ passes through $y$ and $y'$.
\end{enumerate}
Of course the same properties hold in the dual Cayley plane. We can encode
our incidence geometry in the following diagram:
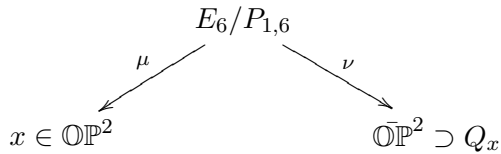
\begin{equation*} 
\xymatrix
{ & E_6/P_{1,6}\ar[dr]^{\nu}\ar[dl]_{\mu} & \\
x\in \OP2 & & \bOP2\supset Q_x }
\end{equation*}
where $P_{1,6}\subset E_6$ is the parabolic subgroup defined as 
the intersection of $P_1$ and $P_6$. In particular $Q_x=\nu(\mu^{-1}(x))$. 

\smallskip
Another useful property that we want to mention is that the Cartan
cubic is a homolo\"idal polynomial, in the sense that its derivatives 
define a birational transformation of the projective space. More 
intrinsically, we have a birational quadratic map
\begin{equation*} 
\begin{array}{rcl}
P  Det : \PP\JO & \dashrightarrow & \PP\JO^\vee \\
  x & \mapsto & Det(x,x,\ast),
 \end{array}
\end{equation*} 
compatible with the $E_6$-action. To be precise, the polarized determinant 
$P  Det$  defines 
an isomorphism between the open orbits on both sides, but it contracts
the Cartan cubic to the dual Cayley plane. Moreover it blows-up the Cayley plane 
itself, which is its indeterminacy locus, to the dual Cartan cubic, in 
such a way that a point $x\in\OP2$ is blown-up into the linear span of
the corresponding quadric $Q_x$. Of course, the picture being symmetric, 
the inverse of  $P  Det$ must be the quadratic map defined by the 
derivatives of the dual Cartan cubic. 

In order to make explicit computations, it is convenient to identify
$\JO$ and $\JO^\vee$ (as vector spaces) through the following scalar 
product:
$$\langle \begin{pmatrix} 
r           & w            & \bar{v}  \\ 
\bar{w} & s           &  u            \\ 
v        & \bar{u} &  t            \\ 
\end{pmatrix}, 
\begin{pmatrix} 
r'            & w'            & \bar{v}'  \\ 
\bar{w}' & s'           &  u'            \\ 
v'            & \bar{u}' &  t'            \\ 
\end{pmatrix}\rangle =  rr'+ss'+tt'-2\mathrm{Re}(\bar{u}u'+\bar{v}v'+\bar{w}w').$$
With this convention, the map $P  Det$ is given explicitly by
\begin{equation}\label{PDet}
x=\begin{pmatrix} 
r            & w            & \bar{v}  \\ 
\bar{w} & s           &  u            \\ 
v            & \bar{u} &  t            \\ 
\end{pmatrix} \quad\mapsto\quad
\begin{pmatrix} 
st-|u|^2            &  tw-\bar{v}\bar{u}            & wu-s\bar{v}  \\ 
t\bar{w}-uv &      rt-|v|^2        &    ru-\bar{w}\bar{v}                     \\ 
\bar{u}\bar{w}-sv          & r\bar{u}-vw &       rs-|w|^2           \\ 
\end{pmatrix} .
\end{equation}
The main interest of this presentation
 is that the equation of the dual 
Cartan cubic is given exactly by the same formula as $Det$. In particular 
the above expression of $P  Det$, which is a kind of {\it comatrix} map, 
 defines a birational involution  of $\PP^{26}$. 

\subsection{Cartan representations}

Given a cubic hypersurface $X$ of dimension $d$, we define 
a {\it Cartan representation} of $X$ as a linear section of
the Cartan cubic isomorphic to $X$:
$$X\simeq \mathcal{C}\cap L \qquad \mathrm{for}\quad \PP^{d+1}\simeq 
L\subset\PP\JO.$$
Note that the Cartan cubic being singular in codimension nine, 
a {\it smooth} cubic hypersurface of dimension bigger than eight 
cannot have any Cartan representation. The main result of this 
section is that, on the contrary,  
a general cubic hypersurface of dimension smaller 
than eight always admits a Cartan representation. (The boundary case 
of eight dimensional sections will be considered later). A more precise 
statement is the following. 

\begin{prop}\label{finiterep}
A generic cubic sevenfold is, up to isomorphism, 
a linear section of the Cartan cubic, in a finite, non zero, 
number of different ways. 
\end{prop}

\proof
Let ${\mathcal M}^7_3$ denote the moduli space of cubic sevenfolds. 
What we need to prove is that the rational map
$$\Psi : G(9,\JO)//E_6\dashrightarrow {\mathcal M}^7_3,$$ 
obtained by sending $L\in G(9,\JO)$ to the isomorphism class of 
$X_L=\mathcal{C}\cap L$,
is generically finite. The two moduli spaces have the same dimension so we just
need to check that $\Psi$ is generically \'etale, hence to find an $L$ at 
which its differential is surjective. This will follow from the surjectivity of the map
$$\varphi\in Hom(L,\JO)\mapsto P_\varphi(x)=Det(x,x,\varphi(x))\in S^3L^\vee .$$

We will choose our $L$ as the image of a map
$$(r,y)\mapsto \begin{pmatrix} r & \alpha(y) & \bar{y} \\ \bar{\alpha}(y) & r & \beta(y) \\ y & 
\bar{\beta}(y) & r \end{pmatrix},$$
where $\alpha,\beta:\mathbb{O}\rightarrow \mathbb{O}$ are linear endomorphisms (of $\mathbb{O}$ considered
as an eight-dimensional vector space). 

Consider the morphism $\varphi$ from $L$ to $\JO$ defined by 
$$(r,y)\mapsto \begin{pmatrix} a & w & \bar{v} \\ \bar{w} & b & u \\ v & \bar{u} & c \end{pmatrix},$$
for some linear forms $a,b,c$ and linear functional $u,v,w$. The corresponding polynomial 
$P_\varphi(x)$ is
\begin{eqnarray*}
P_\varphi(x) & = & r^2(a+b+c)-a|\alpha(y)|^2-b|y|^2-c|\beta(y)|^2 \\
 & & -2r Re(u\beta(y)+vy+w\alpha(y)) +2Re (w\beta(y)y+\beta(y)v\alpha(y)+uy\alpha(y)).
\end{eqnarray*}
So what we need to prove is that the $27$ following quadrics in $r,y$ generate the whole 
space of cubics
\begin{eqnarray*}
 &r^2-|\alpha(y)|^2, & \qquad y\alpha(y)-r\beta(y), \\ 
 &r^2-|y|^2, & \qquad\alpha(y)\beta(y)-ry, \\
 &r^2-|\beta(y)|^2, & \qquad\beta(y)y-r\alpha(y).
\end{eqnarray*}

The following easy lemma is left to the reader. 

\begin{lemm}
A sufficient condition for this to be true is that the sixteen quadrics $y\alpha(y)$ and $\beta(y)y$
generate the  whole space of cubics in $y\in \mathbb{O}$. 
\end{lemm}

\medskip To finish the proof, we thus only need to exhibit $\alpha(y)$ and $\beta(y)$ for which the
lemma does hold. Let us choose them of the form 
\begin{eqnarray*}
 \alpha(y)=\sum_is_iy_{\sigma(i)}e_i, \qquad \beta(y)=\sum_it_iy_{\tau(i)}e_i,
\end{eqnarray*}
for some coefficients $s_i,t_i$ and permutations $\sigma, \tau$. 
More precisely, let $\sigma(i)=i+1$ and $\beta(i)=i+3$ (using the cyclic order on $0,\ldots ,7$),
and $s=t=(1,2,-1,3,1,-1,2,1)$. A computation with {\sc Macaulay2} shows that the lemma does indeed 
hold in that case. 
\qed

\medskip\noindent {\it Question}. What is the number of Cartan representations 
of a generic cubic sevenfold? (Otherwise said, what is the degree of $\Psi$?)
Is it equal to one, in which case we would get a {\it canonical form} for 
cubic sevenfolds? Or bigger, in which case the derived category would have 
many interesting symmetries, as we will see later on? 

\medskip\noindent {\it Remark}. In \cite{im-fcy} we have shown that the Cayley
plane belongs to a series of four homogeneous spaces for which we expect very
similar phenomena. One of these is the Grassmannian $G(2,10)$, and the 
statement analogous to Proposition \ref{finiterep} in that case is the fact 
that a general quintic threefold (the archetypal Calabi-Yau) admits a finite
number of Pfaffian representations. This was proved by Beauville and Schreyer 
\cite{beau} but the precise number of such representations is not known (although
it may be computable, being a Donaldson-Thomas invariant). The two other cases
are spinor varieties, for which the corresponding statements have not been 
checked yet.

\subsection{A special vector bundle}
The next question to ask is: what kind of additional data do we need on a 
cubic sevenfold, in order to define it as a linear section of the Cartan cubic?
 
In the case of quintic threefolds, it is easy to see that Pfaffian representations
are in correspondence with certain special rank two vector bundles, obtained 
as the restriction of the kernel bundle on the regular part of the 
Pfaffian cubic, which parametrizes skew-symmetric forms of corank two.

On the regular part of the Cartan cubic, there is a natural quadric bundle 
$\cQ$ of relative dimension eight, whose fiber at $p$ is just the corresponding 
entry locus $Q_p$. The linear span $\langle Q_p\rangle$ contains $p$, and 
since $Q_p$ is always smooth, the polar to $p$ 
with respect to $Q_p$, is a hyperplane $\PP(\cE^\vee_p)\subset \langle Q_p\rangle$. 
This defines a vector bundle $\cE^\vee$ of rank nine, which is a sub-bundle 
of the trivial vector bundle with fiber $\JO^\vee$. We will denote by $\cK$ 
the quotient bundle, whose rank is eighteen. 

Our next aim is to understand the special properties of these vector bundles, 
when restricted to a cubic sevenfold $X$ defined as a general linear section 
$X=\mathcal{C}\cap L$ of the Cartan cubic. 

\begin{prop}\label{aCM}
Let $\cE_X$ and $\cK_X$ denote the restrictions of $\cE$ and $\cK$ to $X$. 
\begin{enumerate} 
\item $\cE_X$ is arithmetically Cohen-Macaulay;
\item its Hilbert polynomial is $P_{\cE_X}(k)=27\binom{k+7}{7}$;
\item the space of its global sections is $H^0(X,\cE_X)=\JO$;
\item $\cE_X^\vee =\cE_X(-2)$ and $\cK_X^\vee =\cK_X(-1)$;
\item $\chi(End(\cE_X))=0$.
\end{enumerate} 
\end{prop}

\proof 
 The polarization of $Det$ defines over $\PP\JO^\vee$ a map
$$\JO^\vee\otimes\cO_{\PP\JO^\vee}(-1)\stackrel{\delta}{\longrightarrow} \JO\otimes\cO_{\PP\JO^\vee}.$$
This map is invertible outside $\mathcal{C}$, and we claim that its cokernel $\cF$
has constant rank $9$ over the smooth part $\mathcal{C}_0$ of $\mathcal{C}$ (that is, outside the Cayley plane). 
To check this, by homogeneity we just need to check this claim at one point of  
$\mathcal{C}_0$; we choose 
$$p=\begin{pmatrix} 1 & 0 & 0 \\ 0 & 1 & 0 \\ 0 & 0 & 0 \end{pmatrix}.$$
Note that the entry locus of that point is the quadric
$$Q_p=\Bigg\{\begin{pmatrix} r & z & 0 \\ z & s & 0 \\ 0 & 0 & 0 \end{pmatrix}, \quad 
rs-|z|^2=0\Bigg\} .$$
An easy computation shows that the image of $\delta_p$ is generated by the 
linear forms $r+s,t,x,y$. Its dual can then be identified with the kernel 
of this linear form, that is 
$$\cF^\vee_p=\Bigg\{\begin{pmatrix} r & z & 0 \\ z & -r & 0 \\ 0 & 0 & 0 \end{pmatrix}, 
\quad r\in\CC, z\in\mathbb{O}\Bigg\}.$$
But this is clearly the linear space polar to $p$ with respect to the quadric $Q_p$. 
Hence an identification $\cF^\vee\simeq\cE^\vee$ over $\mathcal{C}_0$. 
In particular, we deduce over the linear span $\PP_X$ of $X$ an exact sequence of 
sheaves 
\begin{eqnarray*}
0\ra \JO^\vee\otimes\cO_{\PP_X}(-1)\ra \JO\otimes\cO_{\PP_X}\ra \cE_X\ra 0.
\end{eqnarray*}
This immediately implies that $\cE_X$ is arithmetically Cohen-Macaulay. 
That $\cE_X^\vee =\cE_X(-2)$
follows from the fact that the quadric bundle $Q_p$ is everywhere non degenerate over 
$\mathcal{C}_0$. Similarly we get $\cK_X^\vee =\cK_X(1)$. Also we directly deduce
 the Hilbert polynomial of $\cE_X$. 

Now we can restrict the previous sequence to $X$, which gives 
\begin{equation}\label{End-exacte}
0\ra \cE_X(-3)\ra \JO^\vee\otimes\cO_X(-1)\ra \JO\otimes\cO_X\ra \cE_X\ra 0.
\end{equation}
This exact sequence is self-dual, up to a twist. We can tensor it by $\cE_X^\vee$
to deduce the 
Hilbert polynomial of $End(\cE_X)$. Indeed, this polynomial verifies the two equations
\begin{eqnarray*}
 P_{End(\cE_X)}(k)-P_{End(\cE_X)}(k-3) &= &27(P_{\cE_X}(k-2)-P_{\cE_X}(k-3)), \\
 P_{End(\cE_X)}(k) &= &-P_{End(\cE_X)}(-k-6).
\end{eqnarray*}
The first equation follows from the exact sequence above, and the second one from
Serre duality, since $\omega_X=\cO_X(-6)$. There is a unique solution to these two equations, 
given by the following formula:
$$P_{End(\cE_X)}(k) = 3^5\binom{k+6}{7}-3^4\binom{k+5}{5}+3^3\binom{k+4}{3}-3^2\binom{k+3}{1}.$$
In particular $\chi_{End(\cE_X)}=P_{End(\cE_X)}(0)=0$.\qed  

\medskip\noindent {\it Remark}. The characteristic classes of $\cE_X$ can easily
be computed from the sequence (\ref{End-exacte}), for example the Chern character
$$ch(\cE_X)=27\frac{1-\exp(-h)}{1-\exp(-3h)},$$
where $h$ is the hyperplane class on $X$. 
One easily deduces the Chern classes of the self-dual bundle $\cE_X(-1)$:
$$c(\cE_X(-1))=1-3h^2+9h^4-39h^6.$$

\subsection{Relation with the normal bundle}
Recall that with our conventions, the Cayley plane $\OP2=E_6/P_1$. In 
particular the category of $E_6$-homogeneous vector bundles on $\OP2$ 
is equivalent to the category of $P_1$-modules. But recall that the 
semi-simple part of $P_1$ is a copy of $Spin_{10}$, whose minimal 
non trivial representations are the ten-dimensional vector representation,
and the two sixteen-dimensional half-spin representations. The corresponding
homogeneous vector bundles are, up to twists, the tangent bundle $T$ 
and cotangent bundle $\Omega$ for the latter two, and the normal bundle for the former. 
To fix notations, we will denote by $\cS$ the irreducible homogeneous 
vector bundle on $\OP2$ such that $H^0(\OP2,\cS)=\JO$. 
\medskip

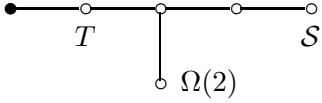
\begin{center}
\setlength{\unitlength}{5mm}
\begin{picture}(20,4)(-5,1)
\multiput(-.3,3.8)(2,0){4}{$\circ$}
\put(-.3,3.8){$\bullet$}
\multiput(0,4)(2,0){4}{\line(1,0){1.7}}
\put(3.7,1.8){$\circ$}
\put(7.7,3.8){$\circ$}
\put(1.6,3){$T$}
\put(4.4,1.8){$\Omega(2)$}
\put(7.6,3){$\cS$}
\put(3.85,2.1){\line(0,1){1.8}}
\end{picture}
\end{center}

The normal bundle is then $\cN=\cS(1)$. 
This implies that the projectivization of the bundle $\cS^\vee$ over $\OP2$ is a desingularization 
of the dual Cayley cubic: the line bundle $\cO_\cS(1)$ is generated by its global sections, whose space 
is isomorphic with $H^0(\OP2,\cS)=\JO$. We have a diagram 
\begin{equation*} 
\xymatrix
{ & \PP(\cS^\vee)\ar[dr]^{\gamma}\ar[dl]_{\pi} & \\
 \OP2 & & \mathcal{C}\supset \bOP2}
\end{equation*}
where the morphism $\gamma$, which is defined by the linear system $|\cO_\cS(1)|$,
  is an isomorphism outside the dual Cayley plane $\bOP2$. (We use the same notation
for the Cartan cubic in the dual space.) Moreover the 
pre-image of $\bOP2$ is the incidence variety $E_6/P_{1,6}$, embedded as a divisor 
in $\PP(\cS^\vee)$. 
Note that the rational map
$\pi\circ\gamma^{-1} : \mathcal{C}\rightarrow \OP2$ is the restriction of the comatrix
map $P Det$. 

Since $X$ does not meet $\bOP2$, $\gamma$ is an isomorphism over $X$, and $\gamma^{-1}(X)$ 
is the complete intersection, in $\PP(\cS^\vee)$, of 18 general sections of $\cO_\cS(1)$.
By construction, we have 
\begin{equation}\label{SvsE}
\pi^*\cS_X^\vee = \cE_X^\vee\oplus \cO_X(-1),
\end{equation}
if we denote by $\pi^*\cS_X$ the restriction of $\pi^*\cS$ to $\gamma^{-1}(X)\simeq X$;
note that the restriction of $\cO_S(1)$ coincides with the pull-back of $\cO_X(1)$ by $\gamma$. 
If we let $H_X$ denote the restriction to $\gamma^{-1}(X)$ of $\pi^*\cO_{\OP2}(1)$, 
we have $H_X\simeq \cO_X(2)$, and the map $\theta_{|X} : X\rightarrow \OP2$ is defined 
by a sub-linear system of $|\cO_X(2)|$.  (One can check that for $X$ generic,  $\theta_{|X}$
is an embedding.)

\begin{theo}\label{rigidity}
The vector bundle $\cE_X$ is simple and infinitesimally rigid. 

More precisely, $$h^i(X,End(\cE_X))=\delta_{i,0}+\delta_{i,3}.$$
\end{theo}
 

\proof 
We start with a simple observation showing that the higher cohomology groups 
of $End(\cE_X)$ cannot all be trivial. 

\begin{lemm}\label{duality}
One has a natural duality $$H^i(X,End(\cE_X))\simeq H^{3-i}(X,End(\cE_X))^\vee .$$
In particular $H^i(X,End(\cE_X))=0$ for $i\ge 4$. 
\end{lemm}

\proof Consider the exact sequence (\ref{End-exacte}). Observe that $\cE_X(-3)$ is acyclic. Indeed, 
 $\cE_X$ is aCM, $H^0(\cE_X^\vee(-3))=0$, and by Serre duality, 
 $H^7(\cE_X(-3))=H^0(\cE_X^\vee(-3))=H^0(\cE_X(-5))=0$. We deduce that  
for any $i$, 
$$H^i(End (\cE_X))=H^{i+2}(End (\cE_X)(-3)).$$ 
But the latter is Serre dual to $H^{5-i}(End (\cE_X)(-3))$ which, 
by the previous assertion, is isomorphic to $H^{3-i}(End (\cE_X))$. \qed

\medskip
Now we use the fact that $\cE_X^\vee=\cE_X(2)$ to decompose $End(\cE_X)$ as the sum of its 
symmetric and skew-symmetric parts $S^2\cE_X^\vee(-2)$ and $\wedge^2\cE_X^\vee(-2)$, which we treat separately. 

\medskip {\bf A.}
We start with the latter. To determine its cohomology, we will use the fact that 
$\wedge^2(\pi^*\cS_X) = \wedge^2\cE_X\oplus \cE_X(1)$.
Moreover, we observe that we can use the Koszul complex 
\begin{equation}\label{koszul}
0\ra \wedge^{18}L\otimes\cO_\cS(-18)\rightarrow\cdots \rightarrow 
L\otimes\cO_\cS(-1)\rightarrow \cO_\cS\rightarrow \cO_X \rightarrow 0
\end{equation}
to deduce the cohomology of $\wedge^2(\pi^*\cS_X)$ from that of $\wedge^2(\pi^*\cS)$ and 
its twists, which we can derive from Bott's theorem. 

Since $\cO_X(-2)=H_X^{-1}$, we will twist the Koszul complex by $\pi^*(\wedge^2\cS\otimes H^{-1})$. 
This bundle is acyclic since $\wedge^2\cS\otimes H^{-1}$ is acyclic on $\OP2$.
Its twists by $\cO_\cS(-k)$ are also acyclic for $k\le 9$, since they are acyclic on the 
fibers of $\pi$. For $k\ge 10$ and any bundle $F$ on $\OP2$, we have an isomorphism
$$\begin{array}{rcl}
H^q(\PP(\cS^\vee),\pi^*F\otimes \cO_\cS(-k)) &= & 
H^{q-9}(\OP2,F\otimes Sym^{k-10}\cS^\vee\otimes \det(\cS^\vee)) \\
&= & H^{q-9}(\OP2,F\otimes Sym^{k-10}\cS\otimes H^{5-k}).
\end{array}$$

We will need a few decomposition formulas. Denote by $\cS_m$ the homogeneous vector bundle 
on $\OP2$ defined by the weight $m\omega_6$. In particular $\cS_1=\cS$. 

\begin{lemm} For any positive integer $m$, 
$$Sym^{m}\cS = \bigoplus_{\ell\ge 0}\cS_{m-2\ell}\otimes H^\ell .$$
\end{lemm}

\proof This follows from the decomposition formulas of the symmetric powers of the 
natural representations of special orthogonal groups, which are classical. \qed 

\medskip Applying the formula above to $F=\wedge^2\cS\otimes H^{-1}$, we thus get 
$$H^q(\PP(\cS^\vee),\pi^*(\wedge^2\cS\otimes H^{-1})\otimes \cO_\cS(-k)) = 
\bigoplus_{\ell\ge 0}H^{q-9}(\OP2,\wedge^2\cS\otimes \cS_{k-2\ell-10}\otimes H^{4-k+\ell}).$$

Then we need to decompose the tensor products $\wedge^2\cS\otimes\cS_m$ into irreducible components.
This decomposition is given by the next Lemma, where  
for future use we include a few more similar formulas. 
We shall  denote by $\cS_{p,q,r}$
the homogeneous vector bundle on $\OP2$ defined by the weight $p\omega_6+q\omega_5+r\omega_4$,
and we let $\cS_{p,q}=\cS_{p,q,0}$.

\begin{lemm}

\begin{enumerate}
\item $\cS_m\otimes \cS = \cS_{m+1}\oplus \cS_{m-1,1}\oplus  \cS_{m-1}\otimes H$.
\item $\cS_m\otimes \wedge^2\cS = \cS_{m,1}\oplus \cS_{m-1,0,1}
\oplus \cS_m\otimes H\oplus \cS_{m-2,1}\otimes H$.
\item $\cS_m\otimes \cS_2 = \cS_{m+2}\oplus \cS_{m,1}\oplus \cS_{m-2,2}
\oplus \cS_m\otimes H\oplus \cS_{m-2,1}\otimes H\oplus \cS_{m-2}\otimes H^2$.
\end{enumerate}
\end{lemm}

\proof Routine representation theory. \qed 

\medskip
Now we can decompose each term $\wedge^2\cS\otimes \cS_{k-2\ell-10}$ into irreducible components,
and compute the cohomology of each term. The special cases of Bott's theorem we will need are covered 
by the following lemma: 

\begin{lemm}\label{coho}
Let $p\ge 0$ and $m\ge 1$. Then
\begin{enumerate}
\item $H^i(\OP2,\cS_p\otimes H^{-m})=0$ for $i<16$, except
$$H^8(\OP2,\cS_p\otimes H^{-m})=V_{(p+4-m)\omega_1+(m-8)\omega_6}\quad \mathrm{{\it if}}\; 8\le m\le p+4.$$
\item $H^i(\OP2,\cS_{p,1}\otimes H^{-m})=0$ for $i<16$, except
$$H^8(\OP2,\cS_{p,1}\otimes H^{-m})=V_{(p+5-m)\omega_1+(m-9)\omega_6}\quad \mathrm{{\it if}}\; 9\le m\le p+5.$$
\item $H^i(\OP2,\cS_{p,2}\otimes H^{-m})=0$ for $i<16$, except
$$\begin{array}{l}
H^4(\OP2,\cS_{p,2}\otimes H^{-5}) = V_{p\omega_1}, \\
H^8(\OP2,\cS_{p,2}\otimes H^{-m})  =  V_{(p+6-m)\omega_1+(m-10)\omega_6}\quad \mathrm{{\it if}}\; 10\le m\le p+6,\\
H^{12}(\OP2,\cS_{p,2}\otimes H^{-p-11}) = V_{p\omega_6}.
\end{array}$$
\item $H^i(\OP2,\cS_{p,0,1}\otimes H^{-m})=0$ for $i<16$, except
$$\begin{array}{l}
H^2(\OP2,\cS_{p,0,1}\otimes H^{-3}) = V_{p\omega_1}, \\
H^6(\OP2,\cS_{p,0,1}\otimes H^{-7}) = V_{(p-2)\omega_1}, \\
H^8(\OP2,\cS_{p,0,1}\otimes H^{-m}) = V_{(p+5-m)\omega_1+(m-10)\omega_6}\quad \mathrm{{\it if}}\; 10\le m\le p+5, \\
H^{10}(\OP2,\cS_{p,0,1}\otimes H^{-p-8}) = V_{(p-2)\omega_6}, \\
H^{14}(\OP2,\cS_{p,0,1}\otimes H^{-p-12}) = V_{p\omega_6}. 
\end{array}$$
\end{enumerate}
\end{lemm}

\proof Apply Bott's theorem as in \cite{derop2}. \qed

\medskip
We deduce the following statement: for $q<25$, 
$$H^q(\PP(\cS^\vee),\pi^*(\wedge^2\cS\otimes H^{-1})\otimes \cO_\cS(-k))
=\delta_{q,19}V_{(k-15)\omega_6}.$$
This implies that for $q<25-18=7$, $H^q(X,\wedge^2\cS_X^\vee(-2))$ is the $q$-th cohomology
group of the complex
$$0\rightarrow  V_{3\omega_6}\rightarrow L\otimes V_{2\omega_6}\rightarrow \wedge^2L\otimes V_{\omega_6}
\rightarrow \wedge^3L\rightarrow 0,$$
going from degree zero to degree five. In particular, $H^0(X,\wedge^2\cS_X^\vee(-2))=0$. 

\begin{prop}
For $L$ generic, the map $L\otimes V_{2\omega_1}\rightarrow V_{3\omega_1}$ is surjective. 

Therefore $H^1(X,\wedge^2\cS_X^\vee(-2))=0$.
\end{prop}

\proof We start with the commutative diagram

$$\begin{array}{ccccc}
 && V_{\omega_1}\otimes V_{\omega_6} && \\
 && \downarrow & & \\
\wedge^2V_{\omega_1}\otimes V_{\omega_1} & \rightarrow & V_{\omega_1}\otimes S^2V_{\omega_1} &
\rightarrow & S^3V_{\omega_1} \\
 || & & \downarrow  &&\downarrow   \\  
\wedge^2V_{\omega_1}\otimes V_{\omega_1} & \rightarrow & V_{\omega_1}\otimes V_{2\omega_1} & \rightarrow & V_{3\omega_1}
\end{array}$$

\medskip
The horizontal complex in the middle is exact, being part of a Koszul complex. 
It surjects to the horizontal complex of the bottom, which is again exact, as one can check from the 
decompositions into irreducible components, as given by LiE \cite{lie}.
The vertical complex in the middle is also exact, the map $V_{\omega_6}=V_{\omega_1}^\vee
\rightarrow S^2V_{\omega_1}$ being given by the differential of the invariant cubic $Det$. 

We deduce from this diagram that the kernel $\bar{K}$ of the map $L\otimes V_{2\omega_1} \rightarrow V_{3\omega_1}$ is 
simply the image of the kernel $K$ of the map
$L\otimes S^2V_{\omega_1} \rightarrow S^3V_{\omega_1}$.   

\begin{lemm}
Let $L\subset V$ be any subspace of dimension $\ell$ and codimension $c$. Then the 
kernel $K$ of the map  $L\otimes S^2V\rightarrow S^3V$ is the image of $\wedge^2L\otimes V$. 
In particular
$$\dim (K)=\frac{\ell(\ell^2-1)}{3}+c\frac{\ell(\ell-1)}{2}.$$
\end{lemm}

\proof Straightforward. \qed

\medskip
In our case, letting $\ell=18$ and $c=9$ we get $\dim (K)=3315$. 

\begin{lemm}
The map $K\rightarrow \bar{K}$ is an isomorphism.
\end{lemm}

\proof We need to check that $K$ does not meet the kernel $L\otimes V_{\omega_6}$ of the projection to 
$L\otimes V_{2\omega_1}$. Otherwise said, we need to check that the map to $Sym^3V_{\omega_1}$ is injective 
on  $L\otimes V_{\omega_6}\subset V_{\omega_1}\otimes S^2V_{\omega_1}$. Let $\ell_i$ be a basis of $L$,
that we complete into a basis of $V_{\omega_1}$ by some vectors $m_j$ generating a supplement $M$ to $L$. 
Let $(y_k)=(\ell_i^\vee, m_j^\vee)$ denote the
dual basis. Any element of $L\otimes V_{\omega_6}$ can be written as $\sum_i\ell_i\otimes x_i$ for some 
vectors $v_i\in V_{\omega_6}=V_{\omega_1}^\vee$. Its image in $L\otimes S^2V_{\omega_1}$ is 
$$\sum_{i,j,k}\ell_i\otimes Det(x_i,y_j,y_k)y_j^\vee y_k^\vee \mapsto \sum_{i,j,k} 
Det(x_i,y_j,y_k)\ell_i y_j^\vee y_k^\vee 
\in S^3V_{\omega_1}.$$
  
\medskip
The decomposition $V_{\omega_1}=L\oplus M$ induces a direct sum decomposition 
$S^3V_{\omega_1}=S^3L\oplus S^2L M\oplus L S^2M\oplus S^3M$. The component of the previous 
tensor on $L S^2M$ is 
$$\sum_{i,j,k} Det(x_i,m_j^\vee,m_k^\vee)\ell_i m_j m_k.$$ 
For it to vanish, we need that $Det(x_i,m_j^\vee,m_k^\vee)=0$ for any $i,j,k$. But this implies that $x_i=0$ 
fo all $i$ since for a generic $L$, one can check that the map 
$$S^2L^{\perp}\hookrightarrow S^2V_{\omega_6}\rightarrow V_{\omega_6}^\vee$$
is surjective. \qed 

\medskip 
We can now complete the proof of the surjectivity of $\phi : L\otimes V_{2\omega_1}\rightarrow V_{3\omega_1}$. 
We know that $\dim V_{2\omega_1}=351$ and $\dim V_{3\omega_1}=3003$.
By the two previous lemmas the kernel of $\phi$ has dimension $3315$, and therefore the dimension of its
image is $18\times 351-3315=3003$. Hence the surjectivity. This completes the proof of the Proposition.\qed 

\medskip {\bf B.} Now we treat $Sym^2\cE_X(-2)$, essentially in the same way. We first observe that 
$Sym^2(\pi^*\cS_X)=Sym^2\cE_X\oplus \pi^*\cS_X(1)$. Recall that $Sym^2(\cS)=\cS_2\oplus H$, 
where the second factor will contribute to the homotheties in $End(\cE_X)$. 
Therefore what we need to compute is the cohomology of $\pi^*(\cS_2\otimes H^{-1})$ restricted to $X$. 
Proceeding as before, we check that     for $q<25$,
$$H^q(\PP(\cS^\vee),\pi^*(        \cS_2\otimes H^{-1})\otimes \cO_\cS(-k))
=\delta_{q,17}V_{(k-12)\omega_6}\oplus \delta_{q,21}V_{(k-18)\omega_6}.$$
This implies that for $q<21-18=3$, the cohomology group $H^q(X,\pi^*(\cS_2\otimes H^{-1})_{|X})$ can be computed as 
the $q$-th cohomology group of the complex
\begin{equation}\label{complex1}
0\rightarrow  V_{6\omega_6}\rightarrow L\otimes V_{5\omega_6}\rightarrow \wedge^2L\otimes V_{4\omega_6}
\rightarrow \wedge^2L\otimes V_{3\omega_6}\rightarrow \wedge^3L,
\end{equation}
where the rightmost term is in degree two. On the other hand, recall that $Sym^2\cE_X(-2)$ is only a
direct factor of  $\pi^*(\cS_2\otimes H^{-1})_{|X}=\pi^*(\cS_2)_{|X}(-2)$, with complementary factor 
$\cS_X(-1)=\cS_X\otimes H_X^{-1}(1)$. After a computation similar to (but easier than) the previous one, we 
get that for $q<25$,
$$H^q(\PP(\cS^\vee),\pi^*(\cS\otimes H^{-1})\otimes \cO_\cS(-1-k))
=\delta_{q,17}V_{(k-12)\omega_6},$$
and we can conclude that for $q<25-18=7$, $H^q(X,\cS_X(-1))$ is  
the $q$-th cohomology group of the complex
\begin{equation}\label{complex2}
0\rightarrow  V_{6\omega_6}\rightarrow L\otimes V_{5\omega_6}\rightarrow \wedge^2L\otimes V_{4\omega_6}
\rightarrow \wedge^2L\otimes V_{3\omega_6}\rightarrow \wedge^3L\rightarrow \cdots ,
\end{equation}
with the same grading as in the complex (\ref{complex1}).
In particular, $H^q(X,\cS_X(-1))$ coincides with $H^q(X,\pi^*(\cS_2\otimes H^{-1})_{|X})$ for $q<3$, and therefore 
$$H^q(X,Sym^2\cE_X^\vee(-2))=\delta_{q,0}\CC \qquad \mathrm{if}\;\; q<3.$$

The proof of the theorem is now complete. Indeed, we have proved that
$$H^0(X,End(\cE_X))= H^0(X,\wedge^2\cE_X(-2))\oplus H^0(X,Sym^2\cE_X(-2))=0\oplus\CC=
 \CC,$$
so that $\cE_X$ is simple, and moreover 
that  $H^1(X,End(\cE_X))=0$. By Lemma \ref{duality} we deduce that $H^2(X,End(\cE_X))=0$
and $H^3(X,End(\cE_X))=\CC$, and we already know that $H^q(X,End(\cE_X))=0$ for $q>3$. \qed

\medskip\noindent {\it Remark}. Observe that the term 
$$H^{21}(\PP(\cS^\vee),\pi^*(\cS_2\otimes H^{-1})\otimes \cO_\cS(-18))=\CC$$
does contribute to $H^3(X,Sym^2\cE_X^\vee(-2))$, in agreement with the fact that  
$H^3(X,End(\cE_X))$ is one-dimensional. 

\medskip\noindent {\it Remark}. It would be interesting to solve the following 
reconstruction problem: being given a smooth cubic sevenfold $X$, and a rank 
nine vector bundle $\cE$ on $X$, having all the properties we have just established
for our bundle $\cE_X$
(and maybe some others), does it necessarily come from a Cartan representation 
of $X$? And if this the case, how can we reconstruct effectively this representation? 
The global sections of $\cE$ provide what should be a copy of $\JO$, and one would
like to reconstruct the cubic determinant just from the bundle. Unfortunately we have
not been able to do that  for $\cE_X$. 

\section{The Tregub-Takeuchi birationality for cubic sevenfolds}\label{tt}

By the results of the previous section, the general cubic sevenfold $X$ has a finite 
number of Cartan representations, where a Cartan representation is the same as 
a presentation of $X$ as a linear section of the Cartan cubic 
$\mathcal{C}\subset \PP\JO = \PP^{26}$ with a subspace $\PP^8$.
In the dual space $\hat{\PP}^{26}=\PP\JO^\vee$, the 18 linear forms defining 
$\PP^8 \subset \PP^{26}$ span a subspace ${\PP}^{17}$ which intersects the 
Cayley plane $\OP2 \subset \hat{\PP}^{26}$ along a sevenfold $Y$, the orthogonal 
linear section of $X$ with respect to its Cartan representation.
In this section we prove that $X$ and $Y$ are birational one
to each other by constructing a birationality which is an analog of the 
Tregub-Takeuchi birationality between the cubic threefold and its general 
orthogonal linear section in $G(2,6)$. 

The analogy is the following.
The general Pfaffian representation of the cubic threefold $X'$ is the same as the
general representation of $X'$ as a linear section of the Pfaffian cubic 
$Pf \subset \PP J_3(\HH)$ with a subspace $\PP^4$, where  
the projective Jordan algebra $\PP J_3(\HH)$ is identified with  
the projective space $\PP(\wedge^2 V)$, where $V = {\CC}^6$.   
The 10 linear forms defining the subspace   $\PP^4 \subset \PP^{14}$ 
span a subspace $\PP^9$ in $\hat{\PP}^{14} = \PP(\wedge^2 V^\vee)$ 
which intersects the quaternionic plane $\HH\PP^2 = G(2,6) \subset \hat{\PP}^{14}$
along a Fano threefold $Y'$ of degree $14$, the orthogonal section of $X'$. 
It is known since G. Fano that $Y'$ are $X'$ are birational. The two types 
of birationalities between $X'$ and $Y'$ -- of Fano-Iskovskikh and 
of Tregub-Takeuchi are related to curves on the cubic threefold $X'$ 
which are hyperplane sections of surfaces with one apparent double point
(OADP), correspondingly del Pezzo surfaces of degree $5$ and rational 
quartic scrolls, see \cite{tr}, \cite{ta}, \cite{cmr}, \cite{ar}.  
Especially the Tregub-Takeuchi birationality,   
or more precisely its inverse $Y' \rightarrow X'$
starts from a point $p \in Y'$, and after a blow-up of $p$ and 
a flop, ends with a contraction of a divisor to a rational normal quartic 
in curve $Y'$.  

Here we describe the analog of the inverse Tregub-Takeuchi birationality 
between the cubic sevenfold $X$ and its dual sevenfold $Y$ defined by 
a Cartan representation of $X$. It is an analog at least because the 
construction described in the proof of Proposition \ref{tto} applied 
to the cubic threefold $X'$ and its dual $Y'$ as above, and after replacing 
$J_3(\OO)$ by $J_3(\HH)$ gives the Tregub-Takeuchi birationality as 
described e.g. in \cite{tr} and \cite{ta} (see the Appendix). 
Just as in the 3-dimensional case, the inverse birationality $Y \rightarrow X$ 
starts with a blow-up of a point $p \in Y$. In the seven dimensional case 
the birationality ends with a contraction of a divisor onto a prime Fano threefold 
$Z \subset X$ of degree $12$, which is a hyperplane section of a fourfold with OADP.
The connection with  varieties with OADP inside the Pfaffian cubic fourfolds 
and Cartan cubic eightfolds is commented in Section \ref{rce}. 
One geometric explanation why the rational quartic curve $C \subset X$ and 
the Fano threefold $Z \subset X$ of degree $12$ appear as indeterminacy loci 
of the birationalities related to points on their dual varieties is that 
the projective tangent space to a point $p$ on $\HH\PP^2 = G(2,6)$ intersects 
on $G(2,6)$ a cone over $\PP^1 \times \PP^3$ (whose general linear section 
of dimension is a rational normal quartic), while the projective 
tangent space to $\OP2$ at $p$ intersects on $\OP2$ a cone over the 
spinor variety $S_{10} = OG(5,10)$ (whose general 3-dimensional linear section 
is a prime Fano threefold of degree $12$), see \cite{lm}.

\subsection{The orthogonal linear sections of the cubic sevenfold on the Cayley plane}
In \cite{im-fcy} we encountered the Cayley plane and its linear sections
as candidates for being Fano manifolds of Calabi-Yau type. Starting from 
a Cartan representation of a general cubic sevenfold, $X=\cC\cap L$, we  
consider the orthogonal section $Y=\OP2\cap L^\perp$. 

\begin{prop}
The variety $Y$ is a Fano manifold of Calabi-Yau type, of dimension 
seven and index three. 
\end{prop}

\proof The Cayley plane has dimension sixteen and index twelve, hence
a general section of codimension nine has dimension seven and index three.
The fact that it is of Calabi-Yau type follows from \cite[Proposition 4.5]{im-fcy}
and Proposition \ref{finiterep}. \qed

\medskip\noindent {\it Remark}. The situation here is very similar to what
happens for cubic fourfolds. Recall that the Pfaffian cubic fourfolds are those that can
be obtained as linear sections of the Pfaffian cubic in $\PP^{14}$. 
The Pfaffian cubic can be defined as the secant variety to $G(2,6)$ (or the 
projective dual variety to the dual Grassmannian $G(4,6)$), and that $G(2,6)$ 
is nothing else than the projective plane $\HP2$ over the quaternions. 

Given such a smooth four dimensional section of the Pfaffian cubic, 
the orthogonal section of the dual Grassmannian is now a smooth K3 surface, 
and the Hodge structure of the cubic fourfold is indeed ``of K3 type''.
The main difference with cubic sevenfolds is that the general cubic fourfold
is not Pfaffian, the Pfaffian ones form a codimension one family in the moduli
space. Nevertheless, the Hodge structure of a  cubic fourfold is always of K3 
type, and this allows in particular to associate to them complete families 
of symplectic fourfolds \cite{bd,ir}.

\subsection{The Tregub-Takeuchi birationalities for the cubic sevenfold}

We will use the basic properties of the incidence geometry on the 
Cayley plane in order to prove that:

\begin{prop}\label{tto}
The general cubic sevenfold $X=\cC\cap L$ and the dual section 
$Y=\OP2\cap L^\perp$ are birationally equivalent. 
\end{prop}

\noindent {\it Remark.} Since $X$ has index six but $Y$ has index three, 
they are not isomorphic. Their Hodge numbers coincide in middle dimension, 
and in fact they do coincide in general because the Cayley plane and the 
projective space have the same Betti (and Hodge) numbers up to degree six.  

\proof We start by fixing a general point $y_0\in Y$. There is a
corresponding  quadric $Q_{y_0}$ on $\bOP2$. 

Let $y\in Y$ be general. The corresponding quadric $Q_y$ meets 
$Q_{y_0}$ along a unique point of $\bOP2$, say $z$. The tangent 
space $T_z$ to the dual Cayley plane at this point is contained in the Cartan
cubic, and we claim that it meets $L$ at a unique point. Indeed,
$T_z$ is a projective space of dimension $16$, and $L$ has codimension 
$18$. But $y$ and $y_0$ belong to $L^\perp$ and they contain 
$T_z$ (recall that $T_z$ is orthogonal to the linear span of 
$\bar{Q}_z$, and that the fact that $z$ belongs to $Q_y$ is equivalent
to the condition that $y$ belongs to $\bar{Q}_z$).
So the intersection of $T_z$ with $L$ is in fact given by
sixteen general conditions, and we end up with a unique point $\psi(y)\in X$. 

Conversely, let $x\in X$ be a general point, then the entry-locus
$Q_x$ meets $Q_{y_0}$ at a unique point of $\bOP2$, that we again call $z$. 
In turn $z$ defines a quadric $\bar{Q}_z$ on the Cayley plane $\OP2$
that we will cut out with $L^{\perp}$. The latter has codimension nine
but it contains $x$, and we claim that the orthogonal hyperplane to $x$
contains $\bar{Q}_z$. In order to check this we have to make a little computation. 
We may suppose that $y_0$ and $PDet(x)=y_1$ are two general points on $\OP2$. 
The group $E_6$ acts transitively on pairs of points in the Cayley plane which 
are not joined by a line contained in $\OP2$, we can therefore suppose that
$$y_0=\begin{pmatrix} 1 & 0 & 0 \\ 0 & 0 & 0 \\ 0 & 0 & 0 \end{pmatrix}, \quad 
y_1=\begin{pmatrix} 0 & 0 & 0 \\ 0 & 0 & 0 \\ 0 & 0 & 1 \end{pmatrix}.$$
This implies that 
$$Q_{y_0}=\begin{pmatrix} 0 & 0 & 0 \\ 0 & * & * \\ 0 & * & * \end{pmatrix}, \quad 
Q_{y_1}=\begin{pmatrix} * & * & 0 \\ * & * & 0 \\ 0 & 0 & 0 \end{pmatrix}.$$
(We have omitted the quadratic conditions expressing that the determinants
of the non-trivial $2\times 2$ blocs in the matrices above, must vanish.)
Hence the intersection point
$$
z=Q_{y_0}\cap Q_{y_1}=\begin{pmatrix} 0 & 0 & 0 \\ 0 & 1 & 0 \\ 0 & 0 & 0 \end{pmatrix},
\quad \mathrm{and} \quad 
Q_z=\begin{pmatrix} * & 0 & * \\ 0 & 0 & 0 \\ * & 0 & * \end{pmatrix}.$$
Finally, recall that a point $x$ such that $PDet(x)=y_1$ must be contained 
in the linear span of $Q_{y_1}$. If we take into account the condition 
that $x$ is orthogonal to $y_0$, we deduce that it must be of the form
$$x=\begin{pmatrix} 0 & * & 0 \\ * & * & 0 \\ 0 & 0 & 0 \end{pmatrix},$$
and therefore $Q_z\subset x^\perp$, as claimed. 

We conclude that the intersection of $\bar{Q}_z$ with $L^\perp$ is in fact given by
eight general linear conditions, and we end up with two points: 
one must be $y_0$, and we denote the other one $\phi(x)$. 

It is then straightforward to check that the rational maps
$\psi$ and $\phi$ between $X$ and $Y$ are inverse one of the other. 
This concludes the proof.  \qed

\medskip\noindent {\it Question.} Note that $\psi$ and $\phi$ factorize
through $Q_{y_0}$, which implies that the general cubic sevenfold
is birational to a complete intersection of bidegree $(2,d)$ in $\PP^9$.
What is $d$?

\medskip Analyzing in detail the structure of the birational maps $\phi$ and $psi$
seems to be rather complicated. The very first steps in this direction are the 
following two statement. 

\begin{prop} 
The rational map from $\OP2$ to $\bOP2$, mapping a general point $y$
to $Q_y\cap Q_{y_0}$, coincides with the double projection $p$ from $y_0$. 
\end{prop}

\proof We may suppose that $y_0$ is the point that if have chosen in the 
proof of Proposition \ref{tto}. Then we can use the local parametrization 
of $\OP2$ around $y_0$ given by \ref{dense}, 
$$y=
\begin{pmatrix} 
1           & w            & \bar{v}  \\ 
\bar{w} & |w|^2            &  \bar{w}\bar{v}          \\ 
v            & vw &  |v|^2               \\ 
\end{pmatrix} \quad \mapsto \quad 
p(y)=
\begin{pmatrix} 
0        & 0            & 0  \\ 
0 & |w|^2            &  \bar{w}\bar{v}    \\ 
0         & vw &  |v|^2               \\ 
\end{pmatrix} .$$
The point $p(y)$ certainly belongs to $Q_{y_0}$. There remains to 
prove that it also belongs to $Q_y$, hence that it is orthogonal to 
the tangent space to the Cayley plane at $y$. This is a
straightforward explicit computation. \qed 

\begin{prop} 
Suppose $X, Y$  and $y_0$ to be general. 
 The indeterminacy locus of $\psi$ is a cone over a smooth canonical
curve of genus seven. The indeterminacy locus of $\phi$ is a smooth prime 
Fano threefold of degree twelve. 
\end{prop}

\proof There are two possible accidents that could prevent $\psi(y)$
from being defined. The first one is that $Q_y$ meets $Q_{y_0}$ in more
than one point. The second one is that $Q_y\cap Q_{y_0}$ is a single 
point $z$ but such that the tangent space $T_z$ meets $L$ in more
than one point. A dimension count shows that generically, the second
situation cannot happen. The first one happens when $y$ belongs to the 
cone $C_{y_0}$ spanned by the lines through $y_0$ in
the Cayley plane. By \cite{lama}, this is a cone over the spinor variety $S_{10}=OG(5,10)$,
inside a sixteen dimensional projectivized half-spin representation. 
Cutting this cone by $L^\perp$, we get a cone over
a codimension nine generic linear section of  $S_{10}$, which is a 
canonical curve of genus seven. 

Let now $x$ be a point in $X$. To prevent $\phi(x)$ from being defined, 
the only possible accident is that $Q_x$ meets $Q_{y_0}$ in more than
one point. If we let $y=PDet(x)$, so that $Q_x=Q_y$, this means that 
$y$ belongs the cone $C_{y_0}$. Let us determine $PDet^{-1}(C_{y_0})$. 
We may suppose that  
$$y_0=\begin{pmatrix} 1 & 0 & 0 \\ 0 & 0 & 0 \\ 0 & 0 & 0 \end{pmatrix}, \quad 
\mathrm{hence} \;\;  C_{y_0}=\Bigg\{
\begin{pmatrix} 1 & a & \bar{b} \\ \bar{a} & 0 & 0 \\ b & 0 & 0 \end{pmatrix},
 \; |a|^2=|b|^2=ab=0
\Bigg\},$$
as we can see from the explicit parametrization of a neighbourhood of $y_0$ 
given in (\ref{dense}). This implies in particular that
the conditions $|a|^2=|b|^2=ab=0$ on the pair $(a,b)\in\OO\oplus\OO$
define a copy of the spinor variety $S_{10}$ in a half-spin representation. 
Now, (\ref{PDet}) shows that $x$ belongs to $PDet^{-1}(C_{y_0})$
if and only if $Det(x)=0$ and 
$$|v|^2=rt, \quad |w|^2=rs, \quad vw=r\bar{u}.$$
If $r\ne 0$, these conditions imply that $wu=s\bar{v}$ and $uv=t\bar{w}$, so  that 
in fact $PDet(x)=y_0$. Generically, this cannot happen in our situation 
since this would imply that $x$ belongs to the linear span of $Q_{y_0}$, 
which is too small to meet $L$. So we must let $r=0$, in which case we 
are left with the conditions $$|v|^2=0, \quad |w|^2=0, \quad vw=0.$$
As we have just seen, this defines a copy of $S_{10}$. Since the parameters
$s,t,u$ remain free, we conclude that $x$ must belong to the join of $S_{10}$
with a $\PP^9$. Cutting this join with $L$ amounts to cut $S_{10}$ along a 
generic three dimensional linear section, which is a prime Fano threefold of 
degree $12$. \qed


\section{Derived categories}\label{derived}

\subsection{The Calabi-Yau subcategory}

Let again $X$ be a smooth cubic sevenfold. Since it is Fano of index $6$, 
the collection $\cO_X,\ldots ,\cO_X(5)$
is exceptional. Denote by $\cA_X$ the full subcategory of $D^b(X)$, 
the derived category of coherent sheaves on $X$, defined as the left semi-orthogonal 
to this exceptional collection. The following statement is a special case of 
\cite[Corollary 4.3]{kuz1}.

\begin{prop}   
$\cA_X$ is a three-dimensional Calabi-Yau category. 
\end{prop}

The terminology {\it non-commutative Calabi-Yau} is sometimes used. 
Kuz\-ne\-tsov has given in \cite{kuz2} interesting examples of non-commutative K3's, 
which are deformations of commutative ones (i.e., of derived categories 
of coherent sheaves on genuine K3 surfaces). Here the situation is a bit
different, since our non-commutative Calabi-Yau cannot be the derived 
category of any Calabi-Yau threefold $Z$, or even a deformation of such a derived
category. Indeed, computing the Hochschild cohomology of $\cA_X$ with the help
of \cite[Corollary 7.5]{kuz3}, we would deduce from the HKR isomorphism theorem
that the Hodge numbers of $Z$ should be such that 
$$\sum_ph^{p,p}(Z)=\sum_ph^{p,p}(X)-6=2 \;\; !$$

\subsection{Spherical twists}

Suppose that $X$ is given with a Cartan representation $X=\cC\cap L$. Recall
that this induces a vector bundle $\cE_X$ of rank nine on $X$, such that $\cE_X(-1)$ is self-dual.   

\begin{prop}
$\cE_X(-1)$ and $\cE_X(-2)$ are two objects of $\cA_X$.
\end{prop} 

\proof Immediate consequence of Proposition \ref{aCM} and Serre duality.\qed
 
\medskip
A nice consequence is that Theorem \ref{rigidity} can now be interpreted as the 
assertion that $\cE_X(-1)$ and $\cE_X(-2)$ are {\it spherical objects} in 
$\cA_X$, in the sense of Seidel and Thomas \cite{ST}. In particular,
with their notations and terminology:

\begin{coro}\label{sphtwists}
The spherical twists $T_{\cE_X(-1)}$ and  $T_{\cE_X(-2)}$ are auto-equi\-va\-lences
of $\cA_X$.
\end{coro} 

These auto-equivalences have infinite order. This can be seen at the level of 
K-theory: the K-theory of $\cA_X$ is free of rank two, hence has a basis given 
by the classes of  $\cE_X(-1)$ and $\cE_X(-2)$. The matrices $\Phi_{\cE_X(-1)}$ 
and  $\Phi_{\cE_X(-2)}$ of the induced automorphisms of $K(\cA_X)$, 
expressed in these basis, are easy to compute. Since $P_{End(\cE_X)}(-1)=9$, 
they are given by 
$$ \Phi_{\cE_X(-1)}=\begin{pmatrix} 1 & 0 \\ -9 & 1 \end{pmatrix}, \qquad
\Phi_{\cE_X(-2)}=\begin{pmatrix} 1 & 9 \\ 0 & 1 \end{pmatrix}.$$

\noindent {\it Remark}. The category $\cA_X$ should have a very interesting 
group of auto-equivalences, all the bigger than the number of Cartan representations 
of the cubic is large. This is of course related to the fact that $\cA_X$ is Calabi-Yau, 
contrary to $D^b(X)$ which, $X$ being Fano, has no interesting auto-equivalence.

\subsection{Homological projective duality}

We expect for cubic sevenfolds certain phenomena that would illustrate the 
principles of {\it homological projective duality} \cite{kuz1}. Starting 
once again from a Cartan representation $X=\cC\cap L$ of a smooth cubic sevenfold, 
we would like to compare the derived category of $X$ with that of the 
orthogonal section $Y=\OP2\cap L^\perp$, that we suppose to be smooth as well. 
 
In order to describe $D^b(Y)$ we will start from the Cayley plane. Recall that 
on $\OP2$ we have denoted by $\cS$ the rank ten spin bundle. It is endowed with
a non degenerate quadratic form 
$$\mathrm{Sym}^2\cS\rightarrow\cO_{\OP2}(1),$$ 
(which define the quadrics on the dual Cayley plane parametrized by $\OP2$), 
whose kernel is an irreducible homogeneous bundle that we denoted $\cS_2$. 
Rephrasing the results of Lemmas 1, 2, 3, 4, 6 in \cite{derop2} 
(which also imply that $\cS$ and $\cS_2$ are exceptional bundles), we 
can  state the following result:

\begin{theo}\label{lefschetz}
In the derived category of coherent sheaves on $\OP2$, let
\begin{eqnarray*}
 &\cA_0 = \cA_1 = \cA_2 = \langle \cS_2^\vee,\cS^\vee,\cO_{\OP2}\rangle, \\
 &\cA_3 = \cdots = \cA_{11}=\langle \cS^\vee,\cO_{\OP2}\rangle .
\end{eqnarray*}
Then $\cA_0, \cA_1(1),\ldots ,\cA_{11}(11)$ is a semi-orthogonal collection in $\cD^b(\OP2)$. 
\end{theo}

This semi-orthogonal collection is certainly complete (otherwise said, it 
should define a {\it Lefschetz decomposition} of $\cD^b(\OP2)$), but we have 
not been able to prove it. 

\smallskip
If we denote by $\cS_Y$ the restriction of $\cS$ to the linear section $Y$, we deduce:
 
\begin{prop}
The sequence $\cO_Y, \cS_Y, \cO_Y(1), \cS_Y(1), \cO_Y(2), \cS_Y(2)$ is exceptional.
\end{prop}

\proof Use the Koszul complex describing $\cO_Y$ and apply Bott's theorem 
on the Cayley plane. More precisely, use Lemma 1 and Lemma 3 in \cite{derop2}.\qed 

\medskip
It is then tempting to consider the right orthogonal $\cA_Y$ in $D^b(Y)$ to this exceptional
collection and ask:
\begin{itemize}
\item Is $\cA_Y$ a three-dimensional Calabi-Yau category? 
\item Is $\cA_Y$ equivalent to $\cA_X$?
\end{itemize}

Note that if $X$ has $d>1$ Cartan representations, the corresponding orthogonal
sections $Y_1,\ldots , Y_d$ would define non-commutative Fourier partners to 
the non-commutative Calabi-Yau threefold $\cA_X$. 
 
\medskip\noindent {\it Remark}. It follows from the results of \cite[Lemma 4]{derop2} that 
$\cS_{2,Y}(-1)$ belongs to $\cA_Y$. Moreover, completing the computations made for proving this Lemma, 
one gets the following statement: if $0\le i\le 13$, 
$$h^q(\OP2,End(\cS_2)(-i))=\delta_{i,0}\delta_{q,0}+\delta_{i,3}\delta_{q,4}+
\delta_{i,6}\delta_{q,8}+\delta_{i,9}\delta_{q,12}.$$
One can easily deduce that the values of $h^q(Y,End(\cS_{2,Y}))$ are equal to $1, 84, 84, 1, 
0, 0, 0, 0$. This is remarkably coherent with the expected Calabi-Yau property of $\cA_Y$.

\section{The rationality of the Cartan cubic eightfolds}\label{rce}

Contrary to cubic sevenfolds, an eight dimensional linear section of the Cartan cubic cannot
be a general cubic eightfold. Indeed, the dimension of this family of cubics is $\dim G(10,27)-\dim E_6=170-78=92$, 
while the moduli space of cubic eightfolds has dimension $\dim |\cO_{\PP^9}(3)|-\dim PGL_{10}=219-99=120$. 

In the spirit of the study in Section \ref{cartan}, 
the Cartan cubic eightfold is the octonionic analog 
of the Pfaffian cubic fourfold, the last corresponding to the quaternions. 
The most important property of the 
Pfaffian cubic fourfold, known since G. Fano is that it is rational, 
see \cite{fa}. In this section we prove that the Cartan cubic 
eightfold is also rational.   

By definition, an $n$-dimensional variety $V \subset {\PP}^{2n+1}$ 
is a {\it variety with one apparent double point} (VOADP) if through 
the general point $x \in \PP^{2n+1}$ passes a unique secant line to $V$. 
For example the union of two non-intersecting
$n$-dimensional linear spaces in $\PP^{2n+1}$ is a  VOADP.

A cubic $2n$-fold $X \subset \PP^{2n+1}$ 
containing a VOADP $V$ is rational. In fact, since through the 
general point of $X$ passes a unique secant line to $V$,
taking the intersection of this secant line with a general 
hyperplane $H$ yields a birationality between $X$ and $H$. 
For example, a smooth cubic surface $S$ is rational 
because it contains space cubic curves, or pairs of non 
intersecting lines, which are examples of VOADP. 
 
By \cite{CG} any smooth cubic threefold is nonrational, 
but it is unknown whether the general cubic hypersurfaces 
of dimension $\ge 4$ are rational or not. 
The singular cubic hypersurfaces in any dimension are 
evidently rational.  
The general Pfaffian cubic 4-fold $X$ is smooth 
and rational since $X$ contains
quintic del Pezzo surfaces and rational quartic scrolls,
which are other examples of VOADP (see e.g. \cite{ar}).  

Other interesting series of examples of smooth rational 
cubic fourfolds have been found by Hassett \cite{Ha},  
but as far as we know, until now the only 
examples of smooth rational cubic 
hypersurfaces of dimension $\ge 5$ are of even 
dimension $2n$ and contain a VOADP. One delicate part in 
these examples is the proof that a general cubic $2n$-fold
through a particular kind of VOADP is indeed smooth.
Moreover the examples of $n$-dimensional VOADP are 
not too many (see \cite{cmr}), which seems to be one more reason why 
there are very few known examples of rational smooth cubic $2n$-folds 
of dimension bigger than four (see e.g. \cite{ar} for a modern treatment 
of Fano's results on Pfaffian cubic fourfolds 
and for examples of smooth cubic sixfolds containing 
a $3$-dimensional VOADP).

We shall see below that the general Cartan cubic $8$-fold $X$ contains 
a $4$-dimensional VOADP $Y$, which implies that $X$ is rational. This special 
VOADP $Y$ is a transversal section of the spinor variety 
$S_{10}$ (Ex. 2.8 in \cite{cmr}).   
By the comments to Section \ref{tt}, the rationality construction 
coming from the degree $12$ VOADP $Y$ on the Cartan cubic 8-fold $X$ 
is the analog of the rationality construction coming from the quartic 
rational normal scroll on the Pfaffian cubic fourfold, as described in 
\cite{fa} and Thm. 4.3 in \cite{ar}.

\begin{theo}\label{rational}
The general Cartan cubic eightfold is rational.
\end{theo}
 
\proof The key observation is the following. 

\begin{lemm}
Any general $\PP^{15}\subset\PP(H)$ is the linear span of a copy 
of the spinor variety $S_{10}$, contained in the Cayley cubic.
\end{lemm}

As we have already mentioned, this implies the theorem as follows. Let $X$ be the cubic eightfold defined as the linear section 
of the Cayley cubic by a general $\PP^9_X\subset \PP(H)$. By the Lemma, there exists a copy $S_{10}$
of the spinor variety, contained in the Cayley cubic, whose linear span contains $\PP^9_X$ as a general
linear subspace. Then the intersection $Y=S_{10}\cap\PP^9_X$ is smooth of dimension four, 
and is a VOADP. 

Now, as we already mentionned, any cubic $X$ of even 
dimension $2m$ containing a $m$-dimensional VOADP $Y$ is rational. Indeed,
consider a general hyperplane $H$ in $\PP^{2m+1}$ and $x$ a general point of $X$. There is a unique 
secant line $\ell_x$ to $Y$ passing through $x$ and we can let $\pi(x)=\ell_x\cap H$. Then $\pi$ is 
a birational isomorphism whose inverse is described in much the same way: if $h$ is a general point 
of $H$, there is a unique secant line $d_h$ to $Y$ passing through $h$; this line $d_h$ cuts $X$ in two 
points of $Y$, and a third point which is $\pi^{-1}(h)$. \qed

\medskip\noindent {\it Proof of the Lemma}. 
Recall that the spinor variety $S_{10}$ parametrizes the space of lines 
in the Cayley plane passing through a given point. We have deduced that 
$$S_{10}\simeq \PP\{u+v\in\OO\oplus\OO, \;\; |u|^2=|v|^2=u\bar{v}=0\}\subset\PP(\OO\oplus\OO)\simeq\PP^{15}.$$
Let $\rho,\tau$ be linear forms on $\OO\oplus\OO$, and let $\alpha,\beta$ be endomorphisms of $\OO$. Then, 
letting $y=\alpha(u)+\tau(v)$, the set of matrices of the form
$$\begin{pmatrix} \rho(u,v) & u & \bar{y} \\ \bar{u} & 0 & \bar{v} \\ y & v & \sigma(u,v)
\end{pmatrix}, \qquad |u|^2=|v|^2=u\bar{v}=0,$$
defines a copy of $S_{10}$ contained in the Cayley cubic. The corresponding family $F_0$ of linear 
spaces is a smooth subset of $G(16,H)$; to be precise, consider
the decomposition $H=L_0\oplus L_1$, where 
$$L_0=\Bigg\{\begin{pmatrix} 0 & u & 0 \\ \bar{u} & 0 & \bar{v} \\ 0 & v & 0
\end{pmatrix}\Bigg\}, \qquad 
L_1=\Bigg\{\begin{pmatrix} r & 0 & y \\ 0 & s & 0 \\ \bar{y} & 0 & t 
\end{pmatrix}\Bigg\}.$$
Then the set of points in $G(16,H)$ defined by subspaces tranverse to $y_1$ form an 
affine neighbourhood of $L_0$ isomorphic to $Hom(L_0,L_1)$, inside which $F_0$ is
the linear subspace $Hom(L_0,L'_1)$, if $L'_1\subset L_1$ denotes the hyperplane
defined by $s=0$. Consider the map
$$\begin{array}{rrcl}
 \psi :  & E_6\times F_0 & \rightarrow  & G(16,H) \\
  & (g,L)& \mapsto & g(L).
\end{array}$$
The image of $\psi$ consists in linear spaces spanned by a copy of the spinor variety 
contained in the Cayley cubic, since this is the case for $F_0$ and this property is 
preserved by the action of $E_6$. 
We claim that $\psi$ is dominant, which implies the Lemma. More precisely, we will
check that the differential of $\psi$ at $(1,L_0)$
is surjective. To see this first observe that $Im(\psi_*)\subset T_{L_0}G(16,H)\simeq Hom(L_0,L_1)$
contains $T_{L_0}F_0\simeq Hom(L_0,L'_1)$. It also contains the image $\psi_*(\fe_6)$ of 
the Lie algebra $\fe_6$, obtained by 
differentiating the restriction of $\psi$ to $E_6\times \{L_0\}$. So we just need
to prove that $Hom(L_0,L'_1)$ and $\psi_*(\fe_6)$ do span $Hom(L_0,L_1)$, which is equivalent 
to the fact that the quotient map $\fe_6\ra Hom(L_0,L_1/L'1)\simeq L_0^\vee$
is surjective. 

To check this, remember that $\fe_6\simeq\ff_4\oplus H_0$, where $H_0\subset H$ is the hyperplane
of traceless matrices, action on $H$ by Jordan multiplication:
$$M\in H_0\mapsto T_M\in End(H), \qquad T_M(X)=\frac{1}{2}(XM+MX).$$
In particular, for a traceless matrix  
$$M=\begin{pmatrix} 
\alpha & c & \bar{b} \\ \bar{c} & \beta & a  \\ b & \bar{a} & \gamma
\end{pmatrix},$$
a short computation yields
$$T_M(L_0)=\Bigg\{\begin{pmatrix} * & * & * \\ * & \langle a,u\rangle +\langle c,v\rangle
 & * \\ * & * & * \end{pmatrix}, \quad u,v\in\OO\Bigg\}.$$
The central coefficient defines the restriction to $H_0\subset\fe_6$ of the map to 
$L_0^\vee$ we are interested in. It is obviously surjective, hence we are done. \qed
   

\section{Appendix}

We have shown in Proposition \ref{tto} that the general cubic sevenfold, considered as 
a linear section of the Cartan cubic, is birational to the orthogonal section 
of the Cayley plane. In this appendix, we show that the same construction,
when octonions are replaced by quaternions, allows to recover the Tregub-Takeuchi 
birationality between a Pfaffian cubic threefold and the orthogonal section 
of the Grassmannian $G(2,6)$. 

First we recall the Tregub-Takeuchi birationality, see \cite{ta}. 
Let $Y = Y_{14}$ be a smooth prime Fano threefold of degree $14$, 
and let $y_0 \in Y$ be a point through which do not pass any of the 
one-dimensional family of lines on $Y$. 
The linear system $|2 - 5y_0|$ 
defines a birational isomorphism $\psi_{T}: Y \rightarrow X$
with a cubic threefold $X$, the {\it Tregub-Takeuchi} birationality. 
It can be decomposed as
$$\xymatrix{
 & Y'\ar@{.>}[rr]\ar[dl]_{\sigma}\ar[dr] & & Y^+\ar[dr]\ar[dl] & \\
Y\ar@{.>}[rr]^{p}  & & \bar{Y} &  &X
}$$
where $p : Y \dashrightarrow \bar{Y}$ is the double projection through $y_0$,
$\sigma : Y' \rightarrow Y$ is the blow-up of $y_0$, $Y' \dashrightarrow Y^+$ 
is a flop over $\bar{Y}$, and  $Y^+ \rightarrow X$ is the contraction of 
an irreducible divisor $N^+$ of $Y^+$ onto a rational normal quartic curve 
$\Gamma$ on $X$ (the blow-up of $\Gamma$). 
There is on $Y$ a $\PP^1$-family of 1-cycles whose general 
member is a rational quintic curve $C_t$ with a double point at $y_0$,  
and the proper image on $Y^+$ of the general $C_t$
is  the general fiber of the contraction 
$N^+ \rightarrow \Gamma$. 

\smallskip
We have already mentioned the fact that $G(2,6)=\HP2$, with its ``plane projective
geometry'' defined by the family of copies of $G(2,4)\simeq \QQ^4$ parame\-tri\-zed by 
the dual Grassmannian $G(4,6)$. As Mukai has shown, the general prime Fano threefold
$Y$ can be realized as a linear section $Y = G(2,6) \cap L^\perp$, with $L^\perp
=\PP(V_{10})\simeq \PP^9$. The intersection of the Pfaffian cubic in the dual projective
space, with $L\simeq \PP^4$, is a smooth cubic threefold, and it follows from the 
results of \cite{imr} that it must coincide with the target $X$ of the Tregub-Takeuchi 
birationality.  

On the other hand, we have the birational map $\psi : Y\dasharrow X$ 
similar to the one we constructed
in the similar situation over the octonions. Recall the construction: take a 
general point $y$ on $Y$, then $Q_y\cap Q_{y_0}$ is a single point $z$, and the 
tangent space $T_z$ to the Grassmannian $G(4,6)$ at $z$ meets $L$ in a unique 
point $\psi(y)$. 

\begin{prop}
The birational maps $\psi$ and $\psi_T$ are the same. 
\end{prop}

\proof
First observe that $y$ and $y_0$ define two planes $P$ and $P_0$ in 
$\CC^6$, and the point $z=Q_y\cap Q_{y_0}$ in $G(4,6)$ is nothing
else that $P+P_0$ when these two planes are transverse. Exactly as
in the octonionic case, the rational map $y\mapsto z$ coincides
with the double projection from $y_0$. 

Now we study the $\PP^1$-family of quintic rational curves $C_t\subset Y$, 
$t \in \PP^1$, passing doubly through $y_0$.   
Let $\PP^8_{y_0} = \PP(P_0 \wedge\CC^6)$ be the projective tangent space
to $G(2,6)$ at $y_0$, and let 
$\PP^3_{y_0} = \PP(U_4) = \PP((P_0 \wedge \CC^6) \cap V_{10})$
be the projective tangent space to $Y$ at $y_0$. 
These quintics are obtained as follows. For any hyperplane $U_5\subset\CC^6$, 
the intersection $C  = Y \cap G(2,U_5)$ is a quintic curve of arithmetic 
genus $p_a(C) = 1$, since $G(2,U_5)$ has degree five and index five. In this 
$\PP^5$-family, those having a double point at $y_0$ are the rational
curves of our $\PP^1$-family $C_t$. They correspond to those 
5-spaces $U_t$ such that the projective tangent spaces 
to $G(2,U_t)$ and to $Y$ at $y_0$ intersect each other 
along a plane (the tangent space to $C_t$ at the node $y_0$), that is
$$
dim (U_4 \cap \wedge^2 U_t) = 3.
$$ 
 
\begin{lemm}
Let $C_t$ be a rational quintic on $Y$ passing doubly through $y_0$. 
Then the birationality $\psi: Y \dashrightarrow X$  contracts $C_t$. 
\end{lemm}

\proof 
Since the quintic curve 
$C_t$ has a double point at $y_0$, the double projection 
$\bar{p}$ sends $C_t$ to a line $\ell_t \subset \bar{Y}$. 
A point $y_s \in C_t$, $y_s \not= y_0$ is mapped by $\bar{p}$ 
to a point $z_s = z(y_s,y_0)$  
on the line $\ell_t \subset \bar{Y} \subset Q_{y_0} \subset G(4,6)$.
We will then obtain $\psi(y_s)$ as the intersection point of 
the tangent space $T_{z_s}$ to $G(4,6)$ at $z_s$, with $L\simeq\PP^4$,
and we must check that this point does not depend on $s$. For  this 
we observe that 
$$\cap_{z\in \ell_t}T_z :=T_{\ell_t}$$
is five dimensional. Moreover, we claim that $T_{\ell_t}$ is orthogonal
to the linear span of $C_t$. This will imply the claim, since $\langle
C_t\rangle ^\perp$ is a $\PP^9$ inside which $L\simeq\PP^4$ and 
$T_{\ell_t}\simeq \PP^5$ will meet, generically, at a unique point $x_t$. 
But then $\{\psi(y_s)\}=T_{z_s}\cap L\subset T_{\ell_t}\cap L=\{x_t\}$,
hence $\psi$ contracts $C_t$ to the point $x_t$. 

There remains to prove the previous claim. It is enough to check that
any $y_s\in C_t$ is orthogonal to $T_{z_s}$. Recall that if $y_s$ 
represents a plane $P_s\subset \CC^6$, transverse to $P_0$, then 
$z_s$ represents the four-space $P_s+P_0$. But then $\hat{T}_{z_s}=
\wedge^3(P_s+P_0)\wedge\CC^6\subset P_s\wedge (\wedge^3\CC^6)$, and 
therefore $\wedge^2P_s\wedge\hat{T}_{z_s}=0$, which is the required
orthogonality condition. \qed

\medskip We can finally conclude the proof of the proposition. Our two birational
maps $\psi$ and $\psi_T$ induce two birational morphisms $\psi^+$ and $\psi_T^+$ 
from $Y^+$ to $X$. These two morphisms contract the divisor $N^+$, 
and the fibers of their restrictions to $N^+$ are the same. 
Since the  relative Picard number is one, $\psi^+$ and $\psi_T^+$ cannot contract 
any other divisor. Then $\psi^+\circ (\psi_T^+)^{-1}$ is an isomorphism between
$X-\psi_T^+(N^+)$ and $X-\psi^+(N^+)$, and extends continuously to $X$. 
Hence it must extend to an automorphism of $X$. \qed

\end{document}